\documentclass[11pt]{article}

\usepackage{amssymb,amsmath,latexsym}

\begin{document}

\newtheorem{thm}{Theorem}[section]
\newtheorem{lemma}[thm]{Lemma}
\newtheorem{defn}[thm]{Definition}
\newtheorem{prop}[thm]{Proposition}
\newtheorem{corollary}[thm]{Corollary}
\newtheorem{remark}[thm]{Remark}
\newtheorem{example}[thm]{Example}

\numberwithin{equation}{section}

\def\ee{\varepsilon}
\def\qed{{\hfill $\Box$ \bigskip}}
\def\MM{{\cal M}}
\def\BB{{\cal B}}
\def\LL{{\cal L}}
\def\FF{{\cal F}}
\def\EE{{\cal E}}
\def\QQ{{\cal Q}}
\def\AA{{\cal A}}
\def\CC{{\cal C}}
\def\R{{\bf R}}
\def\N{{\mathbb N}}
\def\E{{\bf E}}
\def\F{{\bf F}}
\def\H{{\bf H}}
\def\P{{\bf P}}
\def\Q{{\bf Q}}
\def\S{{\bf S}}
\def\J{{\bf J}}
\def\K{{\bf K}}
\def\G{{\bf G}}
\def\F{{\bf F}}
\def\A{{\bf A}}
\def\loc{{\bf loc}}
\def\eps{\varepsilon}
\def\semi{{\bf semi}}
\def\wh{\widehat}
\def\pf{\noindent{\bf Proof.} }
\def\dim{{\rm dim}}

\title{\Large \bf Estimates on Green functions and
Schr\"{o}dinger-type equations for
non-symmetric diffusions with measure-valued drifts}

\author{Panki Kim\\
Department of Mathematics\\
Seoul National University\\
Seoul 151-742, Republic of Korea\\
Email: pkim@snu.ac.kr \smallskip \\
URL: www.math.snu.ac.kr/$\sim$pkim\\
Telephone number: 82-2-880-4077\\
Fax number: 82-2-887-4694 \\
 and
\smallskip \\
Renming Song\thanks{The research of this author is supported in part
by a joint
US-Croatia grant INT 0302167.}\\
Department of Mathematics\\
University of Illinois \\
Urbana, IL 61801, USA\\
Email: rsong@math.uiuc.edu\\
URL:www.math.uiuc.edu/$\sim$rsong\\
 Telephone number: 1-217-244-6604\\
Fax number: 1-217-333-9576  }

\maketitle

\bigskip

\begin{abstract}
In this paper, we establish sharp two-sided estimates for the
Green functions of non-symmetric diffusions with
measure-valued drifts in bounded Lipschitz domains.
As consequences of these estimates,
we get a 3G type theorem and
a conditional gauge theorem for these diffusions in bounded
Lipschitz domains.

Informally the Schr\"{o}dinger-type operators we consider are of the
form $L+\mu\cdot\nabla + \nu$ where $L$ is uniformly elliptic,
$\mu$ is a vector-valued signed measure
belonging to $\K_{d, 1}$ and $\nu$ is a signed measure
belonging to $\K_{d, 2}$. In this paper,
we establish two-sided estimates for the heat kernels of
Schr\"{o}dinger-type operators in bounded
$C^{1,1}$-domains and a scale invariant boundary Harnack principle
for the positive harmonic functions with respect to Schr\"odinger-type
operators in bounded Lipschitz domains.
\end{abstract}

\vspace{.6truein}

\noindent {\bf AMS 2000 Mathematics Subject Classification}:
Primary: 58C60, 60J45; Secondary: 35P15, 60G51, 31C25,
\bigskip

\noindent {\bf Keywords and phrases:}
Brownian motion, diffusion, diffusion process, non-symmetric diffusion, Kato class, measure-valued drift, transition density,
Green function, Lipschitz domain, 3G theorem, Schr\"{o}dinger operator, heat kernel, boundary Harnack principle,  harmonic function

\bigskip

\begin{doublespace}

\section{Introduction}

This paper is a natural continuation of \cite{KS, KS1, KS2}, where
diffusion (Brownian motion) with
measure-valued drift was discussed.
For a vector-valued signed measure $\mu$
belonging to $\K_{d, 1}$, a
diffusion with measure-valued drift $\mu$ is a diffusion
process whose generator can be informally written as
$L+\mu\cdot\nabla$.
In this paper we consider
Schr\"{o}dinger-type operators  $L+\mu\cdot\nabla + \nu$ (see below for
the definition) and discuss
their properties.

In this paper we always assume that $d\ge 3$. First we recall
the definition of the Kato class $\K_{d, \alpha}$
for $\alpha\in (0, 2]$.
For any function $f$ on $\R^d$ and $r>0$, we define
$$
M^ \alpha_f(r)=\sup_{x\in \R^d}\int_{|x-y|\le r}\frac{|f|(y)dy}
{|x-y|^{d- \alpha}}, \quad \quad 0 <  \alpha \le 2.
$$
In this paper, we mean, by a signed measure, the difference of two
nonnegative measures at most one of which can have infinite
total mass.
For any signed measure $\nu$ on $\R^d$, we use
$\nu^+$ and $\nu^-$ to denote its positive and negative parts, and
$|\nu|=\nu^++\nu^-$ its total variation.
For any signed measure $\nu$  on $\R^d$ and any $r>0$, we define
$$
M^\alpha_{\nu}(r)=\sup_{x\in \R^d}\int_{|x-y|\le r}\frac{|\nu|(dy)}
{|x-y|^{d- \alpha}}, \quad  \quad 0 < \alpha  \le 2.
$$

\begin{defn}\label{d:kc} Let $0 <  \alpha \le 2.$
We say that a function $f$ on $\R^d$ belongs to the Kato
class $\K_{d, \alpha}$ if
$\lim_{r\downarrow 0}M^\alpha_{f}(r)=0$.
We say that a signed Radon measure $\nu$ on $\R^d$
belongs to the Kato class $\K_{d, \alpha}$ if
$\lim_{r\downarrow 0}M^\alpha_{\nu}(r)=0$.
We say that a $d$-dimensional vector valued function $V =
(V^1, \cdots , V^d)$ on $\R^d$ belongs to the Kato
class $\K_{d, \alpha}$ if  each $V^i$  belongs to the Kato
class $\K_{d, \alpha}$.
We say that a $d$-dimensional vector valued signed Radon measure
$\mu = (\mu^1, \cdots, \mu^d)$ on $\R^d$
belongs to the Kato class $\K_{d, \alpha}$ if each $\mu^i$
belongs to the Kato
class $\K_{d, \alpha}$.
\end{defn}

Rigorously speaking a function $f$ in $\K_{d, \alpha}$
may not give rise to a signed measure $\nu$ in $\K_{d, \alpha}$
since it may not give rise to a signed measure at all. However, for
the sake of simplicity we use the convention that whenever we write
that a signed measure $\nu$ belongs to $\K_{d, \alpha}$
we are implicitly assuming that we are covering the case of
all the functions
in $\K_{d, \alpha}$ as well.

Throughout this paper we assume that $\mu=(\mu^1, \dots, \mu^d)$ is
fixed with each $\mu^i$ being a signed measure on $\R^d$ belonging
to $\K_{d, 1}$.  We also assume that the operator
$L$ is either $L_1$ or $L_2$ where
$$
L_1:= \frac12\sum_{i,j=1}^{d} \partial_i (a_{ij} \partial_j)
\qquad \mbox{and} \qquad
L_2:=\frac12 \sum_{i,j=1}^{d} a_{ij} \partial_i \partial_j.
$$
with $\A:=(a_{ij})$ being $C^{1}$ and uniformly elliptic.
We do not assume that $a_{ij}$ is symmetric.

Informally, when $a_{ij}$ is symmetric,
a diffusion process $X$ in $\R^d$ with drift $\mu$ is a
diffusion process
in $\R^d$ with generator
$ L+\mu\cdot\nabla$.
When each $\mu^i$ is given by $U^i(x)dx$ for some function
$U^i$, $X$ is a diffusion in $\R^d$
with generator $L +U\cdot\nabla$
and it is a solution to the SDE
$dX_t=dY_t+U(X_t)\cdot dt$
where $Y$ is  a diffusion in $\R^d$
with generator $L$.
For a precise definition of a (non-symmetric) diffusion $X$ with drift $\mu$
in $\K_{d, 1}$, we refer to section 6 in \cite{KS1}
and section 1 in \cite{KS2}.
The existence and uniqueness
of $X$ were established in \cite{BC} (see Remark 6.1 in \cite{BC}).
In this paper, we will always use $X$ to denote the diffusion process
 with drift $\mu$.

In \cite{KS, KS1, KS2}, we have already studied
some potential theoretical properties of the process $X$.
More precisely, we have established two-sided estimates for
the heat kernel of the killed diffusion process $X^D$
and sharp two-sided estimates on
the Green function of $X^D$ when $D$ is a bounded $C^{1, 1}$ domain;
proved a scale invariant boundary Harnack principle for the
positive harmonic functions of $X$ in bounded Lipschitz
domains; and identified the Martin boundary $X^D$ in bounded
Lipschitz domains.

In this paper, we will first establish
sharp two-sided estimates for the Green function of $X^D$
when $D$ is a  bounded Lipschitz domain. As consequences of
these estimates, we get a 3G type theorem and a conditional
gauge theorem for $X$ in bounded Lipschitz domains.
We also establish two-sided estimates for the heat kernels of
Schr\"{o}dinger-type operators in bounded
$C^{1,1}$-domains and a scale invariant boundary Harnack principle
for the positive harmonic functions with respect to Schr\"odinger-type
operators in bounded Lipschitz domains.
The results of this paper will be used in proving the intrinsic
ultracontractivity of the Schr\"{o}dinger semigroup of $X^D$ in \cite{KS5}.

Throughout this paper, for two real numbers $a$ and $b$, we denote
$a \wedge b := \min \{ a, b\}$ and $a \vee b := \max \{ a, b\}$.
The distance between $x$ and $\partial D$ is denote by $\rho_D
(x)$. In this paper we will use the following convention: the
values of the constants $r_i$, $i=1 \cdots 6$,
$C_0$, $C_1$, $M$, $M_i$,  $i=1 \cdots 5$,  and $\eps_1$ will remain
the same throughout this paper, while the values of the constants
$c, c_1, c_2, \cdots$  may change from one
appearance to another. In
this paper, we use ``$:=$" to denote a definition, which is  read
as ``is defined to be".

\section{Green function estimates and 3G theorem}

In this section we will establish sharp two-sided
estimates for the Green function and a 3G theorem for
$X$ in bounded Lipschitz domains.  We will first
establish some preliminary results for the Green function $G_D(x,y)$
of $X^D$. Once we have these results, the proof of
the Green function estimates is similar to the ones in
\cite{Bo1}, \cite{CK} and \cite{H}. The main difference is
that the Green function $G_D(x,y)$ is not (quasi-) symmetric.

For any bounded domain $D$, we use $\tau_D$ to denote the first
exit time of $D$,
i.e., $\tau_D=\inf\{t>0: \, X_t\notin D\}$.
Given a bounded domain $D\subset \R^d$, we define
$X^D_t(\omega)=X_t(\omega)$ if $t< \tau_D(\omega)$
and $X^D_t(\omega)=\partial$ if $t\geq  \tau_D(\omega)$,
where $\partial$ is a cemetery state. The process
$X^D$ is called a killed diffusion with drift $\mu$ in $D$.
Throughout this paper,
we use the convention $f(\partial)=0$.

It is shown in \cite{KS1} that, for any bounded domain $D$, $X^D$ has a
jointly continuous and strictly positive transition density
function $q^D(t,x,y)$ (see Theorem 2.4 in \cite{KS1}).
In \cite{KS1}, we also showed that
there exist positive constants $c_1$ and $c_2$ depending on $D$
via its diameter such that
for any $(t, x, y)\in (0, \infty)\times D\times D$,
\begin{equation}\label{est:6.1}
\,q^{D}(t, x, y)\,\le\, c_1
t^{-\frac{d}{2}}      e^{-\frac{c_2|x-y|^2}{2t}}
\end{equation}
(see Lemma 2.5 in \cite{KS1}).
Let $G_D(x, y)$ be the Green function of $X^D$, i.e.,
$$
G_D(x, y):=\int^{\infty}_0 q^D(t, x, y)dt.
$$
By (\ref{est:6.1}), $G_D(x,y)$ is
finite for $x\neq y$ and
\begin{equation}\label{G_bd}
G_D(x, y) \, \le \,
   \frac{c}{|x-y|^{d-2}}
\end{equation}
for some $c=c
($diam$(D))>0$.

From Theorem 3.7 in \cite{KS1}, we see that there exist constants
$r_1=r_1(d, \mu)>0$ and $c=c(d, \mu) >1$ depending on
$\mu$ only via the rate at which $\max_{1 \le i \le d}M_{\mu^i}(r)$
goes to zero such that
for $r \le r_1$, $z \in \R^d$, $x, y\in B(z,r)$,
\begin{equation}\label{e:Green_B}
c^{-1}\, |x-y|^{-d+2} \le\, G_{B(z,r)}(x,y)\,
\le\, c\, |x-y|^{-d+2}, \quad x,y \in \overline{B(z, 2r/3)}.
\end{equation}

\medskip

\begin{defn}\label{d:har}
Suppose $U$ is an open subset of $\R^d$.

\begin{description}
\item{(1)} A Borel  function
$u$ defined on $U$
is said to be harmonic with respect to $X$ in $U$ if
\begin{equation}\label{AV}
u(x)\,= \,\E_x\left[u(X_{\tau_{B}})\right],
\qquad x\in B,
\end{equation}
for every bounded open set $B$ with $\overline{B}\subset U$;

\item{(2)} A Borel function
$u$ defined on $\overline{U}$
is said to be regular harmonic with respect to $X$ in $U$ if
$u$ is harmonic with respect to $X$ in $U$ and
(\ref{AV}) is true for $B=U$.
\end{description}
\end{defn}

\medskip

Every positive harmonic function in a
bounded domain $D$ is continuous in $D$ (see Proposition 2. 10 in \cite{KS1}).
Moreover, for every
open subset $U$ of $D$, we have
\begin{equation}\label{e:GH}
\E_x[G_D(X_{T_{U}},y)]=G_D(x,y), \qquad (x,y) \in D \times U
\end{equation}
where  $T_U:=\inf \{ t > 0: X_t \in U\}$. In particular,
for every $y \in D$ and $\eps > 0$,
$G_D(\,\cdot\,, y)$
is regular harmonic in $D \setminus B(y, \eps)$
with respect to $X$
(see Theorem 2.9 (1) in \cite{KS1}).

We recall here
the  scale invariant Harnack inequality from \cite{KS}.

\medskip

\begin{thm}\label{HP} (Corollary 5.8 in \cite{KS})
There exist $r_2=r_2(d, \mu)> 0$ and $c=c(d, \mu) >0$
depending on
$\mu$ only via the rate at which $\max_{1 \le i \le d}
M^1_{\mu^i}(r)$ goes to zero
such that for  every positive
harmonic function $f$ for $X$ in $B(x_0, r)$ with  $r \in (0 , r_2)$,
we have
$$
\sup_{y \in B(x_0, r/2)}f(y) \le c \inf _{y \in B(x_0, r/2)}f(y)
$$
\end{thm}

\medskip

Recall that $r_1>0$ is  the constant from (\ref{e:Green_B}).

\medskip

\begin{lemma}\label{G:HP}
For any bounded domain $D$,
there exists $c=c(D, \mu) >0$ such that for  every
$r \in (0 , r_1 \wedge r_2]$ and $B(z, r) \subset D$,  we have
for every $ x \in
D \setminus  \overline{B(z, r)}$
\begin{equation}\label{e:GH1}
\sup_{y \in B(z, r/2)} G_{D}(y, x)
\le c \inf _{y \in B(z, r/2)} G_{D}(y, x)
\end{equation}
and
\begin{equation}\label{e:GH2}
\sup_{y \in B(z, r/2)} G_{D}(x, y)
\le c \inf _{y \in B(z, r/2)} G_{D}(x, y)
\end{equation}
\end{lemma}

\pf
Fix $x \in D \setminus  \overline{B(z, r)}$.
Since $G_D(\,\cdot\,, x)$ is harmonic for $X$ in $B(z,r)$, (\ref{e:GH1})
follows from Theorem  \ref{HP}. So we only need to show (\ref{e:GH2}).

Since $r < r_1$, by (\ref{G_bd}) and (\ref{e:Green_B}),
there exist $c_1=c_1(D)>1$ and $c_2=c_2(d) >1$ such that
for every $y, w \in \overline{B(z, \frac{3r}4)}$
$$
c_2^{-1}\, \frac{1}{|w-y|^{d-2}} \,\le\,   G_{B(z,r)}  (w,y)
\,\le\,  G_{D}  (w,y)  \,\le\, c_1 \,\frac{1}{|w-y|^{d-2}}.
$$
Thus for $ w \in \partial {B(z, \frac{3r}4)}$ and $ y_1, y_2 \in
{B(z, \frac{r}2)}$, we have
\begin{equation}\label{e:pp3}
G_{D}  (w,y_1) \,\le\, c_1 \,\left(\frac{|w-y_2|}{|w-y_1|}\right)^{d-2}
\frac{1}{|w-y_2|^{d-2}} \,\le\, 4^{d-2}\, c_2c_1\, G_{D}  (w,y_2).
\end{equation}
On the other hand, by (\ref{e:GH}), we have
\begin{equation}\label{e:pp4}
G_{D}(x,y) \, = \, \E_x\left[G_D(X_{T_{B(z,\frac{3r}4)}},y)\right],
\quad y \in B(z, \frac{r}{2})
\end{equation}
Since $X_{T_{B(z,\frac{3r}4)}} \in \partial {B(z, \frac{3r}4)}$,
combining (\ref{e:pp3})-(\ref{e:pp4}), we get
$$
 G_{D}(x,y_1)  \,\le\,  4^{d-2}\, c_2c_1\,
 \E_x\left[G_{D}(X_{T_{B(z,\frac{3r}4)}},y_2)\right]
\,=\, 4^{d-2}\, c_2 c_1 \, G_{D}(x,y_2)  , \quad y_1, y_2 \in B(z,\frac{r}2)
$$
In fact, (\ref{e:GH2}) is true for every $x \in D$.
\qed

\medskip

Recall that a bounded domain $D$ is said to be Lipschitz
if there is a localization radius
$R_0>0$  and a constant
$\Lambda_0 >0$
such that
for every $Q\in \partial D$, there is a
Lipschitz  function
$\phi_Q: \R^{d-1}\to \R$ satisfying
$| \phi_Q (x)- \phi_Q (z)| \leq \Lambda_0
|x-z|$, and an orthonormal coordinate
system $CS_Q$ with origin at $Q$ such that
$$
B(Q, R_0)\cap D=B(Q, R_0)\cap \{ y=(y_1, \cdots, y_{d-1}, y_d)=:
(\tilde y, y_d)\mbox{ in } CS_Q: y_d > \phi_Q (\tilde y) \}.
$$
The pair $(R_0, \Lambda_0)$ is called the
characteristics of the Lipschitz domain $D$.

Any bounded Lipschitz domain satisfies $\kappa$-fat property:
there exists $\kappa_0 \in (0,1/2]$ depending on $\Lambda_0$ such
that for each
$Q \in \partial D$ and $r \in (0, R_0)$ (by choosing $R_0$ smaller
if necessary),
$D \cap B(Q,r)$ contains a ball $B(A_r(Q),\kappa_0 r)$.

In this section, we fix a bounded
Lipschitz domain $D$ with its
characteristics $(R_0, \Lambda_0)$ and $\kappa_0$. Without
loss of generality, we may assume that the diameter of $D$ is less than $1$.

We recall here the scale invariant boundary Harnack principle for $ X^D$ in
bounded Lipschitz domains from \cite{KS1}.

\medskip

\begin{thm}\label{BHP2}  (Theorem 4.6 in \cite{KS1})
Suppose $D$ is a bounded Lipschitz domain. Then there exist
constants $M_1, c >1$  and $r_3 >0 $, depending on
$\mu$ only via the rate at which $\max_{1 \le i \le d}M^1_{\mu^i}(r)$
goes to zero such that for every $Q \in \partial D$, $r < r_3$ and any
nonnegative functions
$u$ and $v$ which are harmonic with respect to $X^D$ in $D \cap B(Q, M_1r)$
and vanish continuously on
$\partial D  \cap B(Q, M_1 r)$, we have
\begin{equation}\label{e:BHP0}
\frac{ u(x)}{ v(x)}
\, \le \,c\, \frac{ u(y)}{ v(y)}  \quad  \mbox{ for any }
x,y \in D \cap B(Q, {r}).
\end{equation}
\end{thm}

\medskip

For any $Q\in \partial D$, we define
\begin{eqnarray*}
\Delta_Q(r) &:=&\left\{ y\mbox{ in } CS_Q
:\, \phi_Q (\tilde y)+ 2r
\, > \,y_d\, >\, \phi_Q (\tilde y),\,\,
 |\tilde y| < 2(M_1+1)r \right\},\\
\partial_1 \Delta_Q(r)   &:=& \left\{ y\mbox{ in } CS_Q
:\, \phi_Q (\tilde y)+
2r
 \,\ge\, y_d \,>\, \phi_Q
(\tilde y), \,\,|\tilde y| = 2(M_1+1)r \right\},\\
\quad \partial_2 \Delta_Q(r)
&:=& \left\{ y\mbox{ in } CS_Q
: \,\phi_Q (\tilde y)+  2r
\,   =\, y_d
, \,\,|\tilde y| \le  2(M_1+1)r  \right\},
\end{eqnarray*}
where $CS_Q$ is the coordinate system with origin at $Q$
in the definition of Lipschitz domains
and $\phi_Q$ is the Lipschitz function there.
Let $M_2 :=
2(1+M_1)\sqrt{ 1 +\Lambda_0^2}+2$ and
$
r_4:= M_2^{-1}(R_0 \wedge r_1 \wedge r_2 \wedge r_3).
$
 If $z \in \overline{\Delta_Q(r)}$ with
$r \le r_4$,
then
$$
|Q-z| \,\le\,
|( \tilde z, \phi_Q (\tilde z))-(\tilde{z}, 0)| +2r
\,\le\, 2r(1+M_1)\sqrt{ 1 +\Lambda_0^2}+2r = M_2r\,\le\, M_2 r_4\,\le\, R_0 .
$$
So  $\overline{\Delta_Q(r)} \subset  B(Q, M_2r )\cap D \subset
B(Q, R_0)\cap D  $.

\medskip

\begin{lemma}\label{l:BHPi}
There exists  constant $c >1$
such that for every $Q \in \partial D$, $r < r_4$, and any
nonnegative functions
$u$ and $v$ which are harmonic in $D \setminus B(Q,r)$
and vanish continuously on
$\partial D \setminus B(Q,r)$, we have
\begin{equation}\label{e:l:BHP}
\frac{ u(x)}{ u(y)}
\, \le \,c\, \frac{ v(x)}{ v(y)}  \quad  \mbox{ for any }
x, y \in D \setminus B(Q, M_2r).
\end{equation}
\end{lemma}

\pf
Throughout this proof, we fix a point $Q$
on $\partial D$, $r < r_4$, $\Delta_Q(r)$,
$\partial_1 \Delta_Q(r)$ and $\partial_2 \Delta_Q(r)$.
Fix an $\tilde y_0\in \R^{d-1}$ with
$|\tilde y_0| = 2(M_1+1)r$.
Since $|(\tilde y_0,  \phi_Q (\tilde y_0))| > r$,
$u$ and $v$ are harmonic with respect to $X$ in $D  \cap
B((\tilde y_0,  \phi_Q (\tilde y_0)), 2M_1r)$
and vanish continuously on
$\partial  D \cap B((\tilde y_0,  \phi_Q (\tilde y_0)),2M_1r)$.
Therefore by Theorem \ref{BHP2},
\begin{equation}\label{B_11}
\frac{ u(x)}{ u(y)}
\, \le \,c_1\, \frac{ v(x)}{ v(y)}  \quad  \mbox{ for any }
x, y \in    \partial_1 \Delta_Q(r) \mbox{ with } \tilde x = \tilde y
=\tilde y_0,
\end{equation}
for some constant $c_1>0$.
Since dist$( D \cap B(Q,r),\partial_2 \Delta_Q(r)) >c r$ for some
$c:=c(\Lambda_0)$, the Harnack inequality
(Theorem \ref{HP}) and
a Harnack chain argument imply that
there exists a constant $c_2 >1$ such that
\begin{equation}\label{B_12}
c_2^{-1}\,<\,\frac{u(x)}{u(y)},\,\,\frac{v(x)}{v(y)} \,<\,c_2,
\quad  \mbox{ for any } x, y\in \partial_2 \Delta_Q(r).
\end{equation}
In particular, (\ref{B_12}) is true with $y
:=( \tilde y_0, \phi_Q(\tilde y_0)+2r)$, which is also in
$\partial_1 \Delta_Q(r)$. Thus (\ref{B_11}) and (\ref{B_12})
 imply that
\begin{equation}\label{B_14}
c_3^{-1} \frac{ u(x)}{ u(y)}
\, \le \, \frac{ v(x)}{ v(y)}  \le \,c_3\,\frac{ u(x)}{ u(y)},
\quad x,y \in  \partial_1 \Delta_Q(r)
 \cup  \partial_2 \Delta_Q(r)
\end{equation}
for some constant $c_3>0$.
Now, by applying the maximum principle (Lemma 7.2 in \cite{KS}) twice,
we get that (\ref{B_14}) is true for every $x \in  D\setminus\Delta_Q(r)
\supset D \setminus B(Q, M_2r)$.
\qed

\medskip

Combining Theorem \ref{BHP2} and Lemma \ref{l:BHPi}, we get
a uniform boundary Harnack principle for $G_D(x, y)$ in
both variables.
Recall $\kappa_0$ is the $\kappa$-fat constant of $D$.

\medskip

\begin{lemma}\label{l:Green_L}
There exist  constants $c >1$, $M > 1/\kappa_0$ and $r_0 \le r_4$
such that for every $Q \in \partial D$, $r < r_0$,  we have for
$x, y \in D \setminus B(Q, r)$ and $z_1, z_2
\in D \cap B(Q, r/M)$
\begin{equation}\label{e:CG_1}
\frac{ G_{D}(x,z_1)}{ G_{D} (y,z_1)}
\, \le \,c\, \frac{ G_{D} (x,z_2) }{ G_{D} (y,z_2)}
\quad
\mbox{and}
\quad
\frac{ G_{D}(z_1,x)}{ G_{D} (z_1,y)}
\, \le \,c\, \frac{ G_{D} (z_2,x) }{ G_{D} (z_2,y)}.
\end{equation}
\end{lemma}

\medskip

Fix $z_0 \in D$ with $r_0/M < \rho_D(z_0) < r_0$ and let $\eps_1:=
r_0 /(12M)$. For $x,y \in D$, we let $r(x,y): = \rho_D(x) \vee
\rho_D(y)\vee |x-y|$ and
$$
\BB(x,y):=\{ A \in D:\, \rho_D(A) > \frac1{M}r(x,y), \,  |x-A|\vee
|y-A| < 5 r(x,y)  \}
$$
if $r(x,y) <\eps_1 $, and $\BB(x,y):=\{z_0 \}$ otherwise.

By a Harnack chain argument we get the following from (\ref{G_bd}) and
(\ref{e:Green_B}).

\medskip

\begin{lemma}\label{G:Poly}
There exists a positive constant $C_0$ such that $G_D(x,y) \le C_0
|x-y|^{-d+2}$,  for all $x,y \in D$,
and $G_D(x,y) \ge C_0^{-1} |x-y|^{-d+2}$ if
$2 |x-y| \le \rho_D (x) \vee \rho_D (y) $  .
\end{lemma}

\medskip

Let $C_1:=C_0 2^{d-2}  \rho_D(z_0)^{2-d}$. The above lemma implies that
$G_D(\cdot, z_0)$ and $G_D(z_0, \cdot)$ are bounded above by $C_1$ on
$D \setminus B(z_0, \rho_D(z_0)/2)$.
Now we define
$$
g_1(x ):=  G_D(x, z_0) \wedge C_1
\quad \mbox{and} \quad
g_2(y ):=  G_D(z_0, y) \wedge C_1.
$$

Using Lemma \ref{G:HP} and a Harnack chain argument,
we get the following.

\medskip

\begin{lemma}\label{G:HI}
For every $y \in D$ and $x_1, x_2 \in D\setminus B(y,
\rho_D(y)/2)$ with $|x_1 -x_2|
\le k (\rho_D(x_1) \wedge
\rho_D(x_2))$, there exists $c:=c(D, k)$ independent of
$y$ and $x_1, x_2$ such that
\begin{equation}\label{e:G:HI}
G_D(x_1,y) \,\le\, c\, G_D(x_2,y)
\quad
\mbox{and}
\quad
G_D(y,x_1) \,\le\, c\, G_D(y,x_2).
\end{equation}
\end{lemma}

\medskip

The next two lemmas follow easily from the result above.

\medskip

\begin{lemma}\label{G:g1}
There exists $c=c(D)>0$ such that for every $x, y \in D$,
$$
c^{-1}\, g_1(A_1)\, \le\, g_1(A_2) \,\le\, c\, g_1(A_1)
\quad
\mbox{and}
\quad
c^{-1}\, g_2(A_1) \,\le\, g_2(A_2)\, \le\, c\, g_2(A_1), \quad A_1, A_2 \in \BB(x,y).
$$
\end{lemma}

\medskip

\begin{lemma}\label{G:g2}
There exists $c=c(D)>0$ such that for every
$x \in \{ y \in D; \rho_D(y) \ge \eps_1/(8M^3)\} $,
$
c^{-1}  \,\le\, g_i(x) \,\le\, c
$, $i=1,2$.
\end{lemma}

\medskip

Using Lemma \ref{G:HP}, the proof of the next lemma is routine
(for example, see Lemma 6.7 in \cite{CZ}).
So we omit the proof.

\medskip

\begin{lemma}\label{G:g3}
For any given $c_1>0$,
 there exists $c_2=c_2(D, c_1, \mu)>0$ such that for every
$|x-y| \le c_1(\rho_D(x) \wedge \rho_D(y))
$,
$$
 G_D(x,y) \,\ge\, c_2 |x-y|^{-d+2}.
$$
In particular, there exists $c=c(D, \mu)>0$ such that for every
$|x-y| \le (8M^3/\eps_1) (\rho_D(x) \wedge \rho_D(y))
$,
$$
c^{-1}\, |x-y|^{-d+2} \,\le\, G_D(x,y) \,\le\, c\, |x-y|^{-d+2}.
$$
\end{lemma}

\medskip

With the preparations above, the following two-sided estimates for
$G_D$ is
a direct generalization of the estimates of the Green function
for symmetric processes
(see \cite{CK} for a symmetric jump process case).

\medskip

\begin{thm}\label{t:Gest}
There exists $c:=c(D)>0$ such that for every $x, y \in D$
\begin{equation}\label{e:Gest}
c^{-1}\,\frac{g_1(x) g_2(y)}{g_1(A)g_2(A)} \,|x-y|^{-d+2}
\,\le\,G_D(x,y)\,\le\,
c\,\frac{g_1(x) g_2(y)}{g_1(A)g_2(A)}\, |x-y|^{-d+2}
\end{equation}
for every $A \in \BB(x,y)$.
\end{thm}

\pf
Since the proof is an adaptation of the proofs of
Proposition 6  in \cite{Bo1} and Theorem 2.4 in \cite{H},
we only give a sketch of the
proof for the case $\rho_D(x) \le \rho_D(y)\le \frac1{2M} |x-y|$.

In this case, we have $r(x,y)=|x-y|$.
Let $r:= \frac12 (|x-y| \wedge \eps_1)$.
Choose $Q_x, Q_y \in \partial D$ with
$|Q_x-x| =\rho_D(x)$ and $|Q_y-y| =\rho_D(y)$. Pick points
$x_1=A_{r/M}(Q_x)$ and
$y_1=A_{r/M}(Q_y)$ so that $x, x_1 \in B(Q_x, r/M)$
and $y, y_1 \in B(Q_y, r/M)$.
Then one can easily check that
$|z_0-Q_x| \ge r$ and $|y-Q_x| \ge r$. So by the first inequality
in  (\ref{e:CG_1}), we have
$$
c_1^{-1}\,  \frac{G_D(x_1,y)}{g_1(x_1)}\,\le\,
\frac{G_D(x,y)}{g_1(x)} \,\le\, c_1 \frac{G_D(x_1,y)}{g_1(x_1)},
$$
for some $c_1>1$.
On the other hand, since  $|z_0-Q_y| \ge r$ and $|x_1-Q_y| \ge r$,
applying the second inequality in  (\ref{e:CG_1}),
$$
c_1^{-1}\,\frac{G_D(x_1,y_1)}{g_2(y_1)} \,\le\,
\frac{G_D(x_1,y)}{g_2(y)} \,\le\, c_1 \frac{G_D(x_1,y_1)}{g_2(y_1)}.
$$
Putting the four inequalities above
together we get
$$
c_1^{-2}\,\frac{G_D(x_1,y_1)}{g_1(x_1)g_2(y_1)}  \,\le\,
\frac{G_D(x,y)}{g_1(x)g_2(y)} \,\le\, c_1^2\frac{G_D(x_1,y_1)}
{g_1(x_1)g_2(y_1)}.
$$
Moreover, $\frac13|x-y| < |x_1-y_1| < 2 |x-y| $ and
$|x_1 -y_1| \le (8M^3/\eps_1) (\rho_D(x_1) \wedge \rho_D(y_1))$.
Thus by Lemma \ref{G:g3}, we have
$$
\frac{1}{2^{d-2}c_2 c_1^{2}} \,\frac{|x-y|^{-d+2}}
{g_1(x_1)g_2(y_1)}  \,\le\,
 \frac{G_D(x,y)}{g_1(x)g_2(y)} \,\le\,3^{d-2}c_2 c_1^2
\frac{|x-y|^{-d+2}}{g_1(x_1)g_2(y_1)},
$$
for some $c_2>1$.

If $r=\eps_1/2$, then $r(x,y)=|x-y| \ge \eps_1$. Thus
$g_1(A)=g_2(A)=g_1(z_0)=g_2(z_0)=C_1$ and
$\rho_D(x_1), \rho_D(y_1) \ge r/M = \eps_1/(2M)$.
So by Lemma \ref{G:g2},
$$
C_1^{-2} c_3^{-2} \,\le\,\frac{g_1(A)g_2(A)}
{g_1(x_1)g_2(y_1)} \,\le\,C_1^2 c_3^2,
$$
for some $c_3>1$.

If $r<\eps_1/2$, then $r(x,y)=|x-y| < \eps_1$ and
$r=\frac12 r(x,y)$. Hence
$\rho_D(x_1), \rho_D(y_1) \ge r/M = r(x,y) /(2M)$.
Moreover, $|x_1-A|, |y_1-A| \ge 6 r(x,y)$.
So by applying the first inequality in  (\ref{e:G:HI}) to $g_1$, and
the second inequality in  (\ref{e:G:HI}) to $g_2$ (with $k=12M$),
$$
c^{-1}_4 \,\le\,\frac{g_1(A)}{g_1(x_1)} \,\le\,c_4
\quad
\mbox{and}
\quad
c^{-1}_4
 \,\le\,\frac{g_2(A)}{g_2(y_1)} \,\le\,c_4
$$
for some constant $c_4=c_4(D)>0$.
\qed

\medskip

\begin{lemma}\label{C:c_L} (Carleson's estimate)
For any given $0 <N< 1$, there exists  constant $c >1$
 such that for every $Q \in \partial D$,  $r < r_0$,
$x \in D  \setminus B(Q, r)$ and $z_1, z_2
\in  D \cap B(Q, r/M)$
with
$B(z_2, Nr) \subset D \cap B(Q, r/M)$
\begin{equation}\label{e:CG_3}
G_{D}(x,z_1)
\, \le \,c\,  G_{D} (x,z_2)
\quad
\mbox{and}\quad
 G_{D}(z_1,x)
\, \le \,c\,  G_{D} (z_2,x)
\end{equation}
\end{lemma}

\pf
Recall that $CS_Q$ is the  coordinate
system with origin at $Q$ in the definition of Lipschitz domains.
Let $\overline{y}:=(\tilde 0, r)$. Since $z_1, z_2
\in  D \cap B(Q, r/M)$, by (\ref{G_bd}),
$$
G_{D} (\overline{y},z_1)
\,\le\, c_1\, r^{-d+2}
\quad
\mbox{and}
\quad
 G_{D} (z_1,\overline{y}) \,
\le\, c_1\, r^{-d+2},
$$
for some constant $c_1>0$.
On the other hand, since $\rho_D(\overline{y}) \ge c_2 r$
for some constant $c_2>0$
and $\rho_D(z_2) \ge  N r$,  by Lemma \ref{G:g3},
$$ G_{D} (\overline{y},z_2)
\,\ge\, c_3\, |\overline{y}-z_2|^{-d+2}
\,\ge \,c_4\, r^{-d+2}
 \quad
\mbox{and}
\quad\, G_{D} (z_2,\overline{y})
\,\ge\, c_3\, |\overline{y}-z_2|^{-d+2}
\,\ge \,c_4\, r^{-d+2},
$$
for some constants $c_3, c_4>0$.
Thus from (\ref{e:CG_1}) with $y=\overline{y}$, we get
$$
G_{D}(x,z_1) \,\le\, c_5\,
 \left(\frac{c_1}{c_4}\right) G_{D} (x,z_2)
\quad
\mbox{and}\quad
G_{D}(z_1,x) \,\le\, c_5\,
\left(\frac{c_1}{c_4}\right) G_{D} (z_2,x)
$$
for some constant $c_5>0$.
\qed

\medskip

Recall that, for $r \in (0, R_0)$,  $A_{r}(Q)$ is
a point in $D \cap B(Q,r)$ such that
$B(A_r(Q),\kappa_0 r) \subset D \cap B(Q,r)$.
For every $x,y \in D$, we denote $Q_x$, $Q_y$ by points on
$\partial D$ such that $\rho_D(x)=|x-Q_x|$ and
$\rho_D(y)=|y-Q_y|$ respectively. It is easy to check that  if $r(x,y) <
\eps_1$
\begin{equation}\label{e:AinB}
A_{r(x,y)}(Q_x),\, A_{r(x,y)}(Q_y) \,\in\, \BB(x,y).
\end{equation}
In fact, by the definition of $ A_{r(x,y)}(Q_x)$,
$\rho_D(A_{r(x,y)}(Q_x)) \ge \kappa_0 r(x,y) > r(x,y) /M$.
Moreover,
$$
|x- A_{r(x,y)}(Q_x)|\, \le\, | x-Q_x|+|Q_x -  A_{r(x,y)}(Q_x)|
\,\le\, \rho_D(x) + r(x,y) \,\le\, 2 r(x,y)
$$
and $ |y- A_{r(x,y)}(Q_x)| \le |x-y| +|x- A_{r(x,y)}(Q_x)| \le 3
r(x,y). $

\medskip

\begin{lemma}\label{G:g4}
There exists $c>0$ such that the following holds:
\begin{description}
\item{(1)}
If $Q \in \partial D$, $0 <s \le r < \eps_1$ and $A=A_r(Q)$, then
$$
g_i(x) \,\le\, c\, g_i(A) \quad \mbox{for every } x \in D
\cap B(Q, Ms) \cap
\{ y \in D : \rho_D(y) > \frac{s}{M} \},\quad  i=1,2.
$$
\item{(2)}
If $x,y,z \in D$ satisfy $ |x-z| \le |y-z|$, then
$$
g_i(A) \,\le\, c\, g_i(B) \quad \mbox{for every } (A, B) \in
\BB(x,y) \times \BB(y,z),\quad  i=1,2.
$$
\end{description}
\end{lemma}

\pf
This is an easy consequence of the Carleson's estimates (Lemma \ref{C:c_L}),
(\ref{e:AinB}) and
Lemmas \ref{G:g1}-\ref{G:g3} (see
page 467 in \cite{H}). Since the proof is similar to the proof
on page 467 in \cite{H},
 we omit the details
\qed

\medskip

The next result is called a generalized triangle property.

\begin{thm}\label{t:tri}
There exists a constant $c >0$ such that for every  $x, y, z \in D$,
\begin{equation}\label{e:tri}
 \frac{G_D(x,y) G_{D}(y,z)}
{ G_{D}(x,z)} \,\le\, c\, \left( \frac{g_1(y)}{g_1(x)}
G_D(x,y) \vee \frac{g_2(y)}{g_2(z)} G_D(y,z)    \right)
\end{equation}
\end{thm}

\pf Let $A_{x,y} \in \BB(x,y)$, $A_{y,z} \in \BB(y,z)$ and
$A_{z,x} \in \BB(z,x)$.
If $|x-y| \le |y-z|$ then $|x-z| \le |x-y|+|y-z| \le 2 |y-z|$.
So by (\ref{e:Gest}) and Lemma \ref{G:g4} (2),
we have
$$
\frac{ G_{D}(y,z)}
{ G_{D}(x,z)} \le c_1^2\frac{g_1(A_{x,z})g_2(A_{x,z}) }
{g_1(A_{y,z})g_2(A_{y,z})}\, \frac{|x-z|^{d-2}}{|y-z|^{d-2}}
\frac{g_1(y)}{g_1(x)} \,\le\, c_1^2 c_2 2^{d-2} \, \frac{g_1(y)}{g_1(x)}
$$
for some $c_1, c_2>0$.
Similarly if $|x-y| \ge |y-z|$, then
$$
\frac{ G_{D}(x,y)}
{ G_{D}(x,z)} \le c_1^2\frac{g_1(A_{x,z})g_2(A_{x,z}) }
{g_1(A_{x,y})g_2(A_{x,y})}\, \frac{|x-z|^{d-2}}{|x-y|^{d-2}}
\frac{g_1(y)}{g_1(x)} \,\le\, c_1^2 c_2 2^{d-2} \, \frac{g_2(y)}{g_2(z)}.
$$
Thus
$$
\frac{G_D(x,y) G_{D}(y,z)}
{ G_{D}(x,z)} \,\le\,c_1^2 c_2 2^{d-2} \, \left( \frac{g_1(y)}
{g_1(x)} G_D(x,y) \vee \frac{g_2(y)}{g_2(z)} G_D(y,z)    \right).
$$
\qed

\medskip

\begin{lemma}\label{G:g5}
There exists $c>0$ such that
for every $x,y \in D$ and $A \in \BB(x,y)$,
$$
g_i(x) \,\vee\, g_i(y)  \,\le\, c\, g_i(A),
 \quad  i=1,2.
$$
\end{lemma}

\pf
If $r(x,y) \ge \eps_1$, the lemma is clear.
If
$r(x,y) <\eps_1$, from Lemma \ref{G:g4} (1), it is easy to see that
that
$$ g_i(x) \, \le\, c\, g_i (A_{r(x,y)}(Q_x))
$$
for some $c>0$, where $Q_x$ is a point on
$\partial D$ such that $\rho_D(x)=|x-Q_x|$.
Thus the lemma follows from
Lemmas \ref{G:g1} and
and (\ref{e:AinB}).
\qed

\medskip

Now we are ready to prove the 3G theorem.

\begin{thm}\label{t:3G}
There exists a constant $c >0$
 such that for every  $x, y, z \in D$,
\begin{equation}\label{3G_est}
 \frac{G_D(x,y) G_{D}(y,z)}
{ G_{D}(x,z)} \,\le\, c\, \frac{|x-z|^{d-2}}{|x-y|^{d-2} |y-z|^{d-2}}.
\end{equation}
\end{thm}

\pf
Let $A_{x,y} \in \BB(x,y)$, $A_{y,z} \in \BB(y,z)$ and $A_{z,x} \in \BB(z,x)$.
By (\ref{e:Gest}), the left-hand side of (\ref{3G_est}) is less
than and equal to
$$
\left(\frac{g_1(y)g_1(A_{x,z})}{g_1(A_{x,y})  g_1(A_{y,z})}\right)
\left(\frac{g_2(y)g_2(A_{x,z})}{g_2(A_{x,y})  g_2(A_{y,z})}\right)
\frac{|x-z|^{d-2}}{|x-y|^{d-2} |y-z|^{d-2}}.
$$
If $|x-y| \le |y-z|$, by Lemma \ref{G:g4} and Lemma \ref{G:g5}, we have
$$
\frac{g_1(y)}{g_1(A_{x,y}) } \le\, c_1,\quad
 \frac{g_2(y)}{g_2(A_{x,y}) } \le\, c_1,\quad
\frac{g_1(A_{x,z})}{  g_1(A_{y,z})}\le\, c_2
\quad
\mbox{and}
\quad
\frac{g_2(A_{x,z})}
{g_2(A_{y,z})}\le\, c_2
$$
for some constants $c_1, c_2>0$.
Similarly, if $|x-y| \ge |y-z|$, then
$$
\frac{g_1(y)}{  g_1(A_{y,z})}\le\, c_1,\quad
\frac{g_2(y)}{  g_2(A_{y,z})} \le\, c_1,\quad
\frac{g_1(A_{x,z})}{g_1(A_{x,y}) }\le\, c_2
\quad
\mbox{and}
\quad
\frac{g_2(A_{x,z})}
{g_2(A_{x,y}) }\le\, c_2.
$$
\qed

\medskip

Combining the main results of this section, we get the following inequality.

\medskip

\begin{thm}\label{t:3G2}
There exist constants $c_1, c_2 >0$
 such that for every  $x, y, z \in D$,
\begin{equation}\label{3G_est2}
\frac{G_D(x,y) G_{D}(y,z)}
{ G_{D}(x,z)} \, \le\,
c_1 \left( \frac{g_1(y)}{g_1(x)} G_D(x,y) \vee
\frac{g_2(y)}{g_2(z)} G_D(y,z)    \right) \le\, c_2\,
\left( |x-y|^{-d+2} \vee  |y-z|^{-d+2}  \right).
\end{equation}
\end{thm}

\pf
We only need to prove the second inequality. Applying Theorem \ref{t:Gest},
we get that there exists $c_1>0$ such that
$$
\frac{g_1(y)}{g_1(x)} G_D(x,y) \,\le\, c_1\,  \frac{g_1(y)g_2(y)}
{g_1(A)g_2(A)} |x-y|^{-d+2}
$$
and
$$
\frac{g_2(y)}{g_2(z)} G_D(y,z) \,\le\, c_1\,  \frac{g_1(y)g_2(y)}
{g_1(B)g_2(B)} |x-y|^{-d+2}
$$
for every $(A,B) \in \BB(x,y) \times \BB(y,z)$. Applying
Lemma \ref{G:g5}, we arrive at the desired assertion.

\section{Schr\"{o}dinger semigroups for $X^D$}

In this section, we will assume that $D$ is a bounded Lipschitz domain.
We first recall some notions from \cite{KS2}.
A measure $\nu$ on $D$ is said to be a smooth measure
of $X^D$ if there is a positive continuous additive functional
(PCAF in abbreviation) $A$ of $X^D$ such that for any
$x\in D$,  $t>0$ and
bounded nonnegative function $f$ on $D$,
\begin{equation}\label{eqn:revuz}
\E_x \int^t_0f(X^D_s) dA_s=\int^t_0\int_D
q^D(s, x, y)f(y) \nu (dy)ds.
\end{equation}
The additive functional $A$ is called the PCAF of
$X^D$ with Revuz measure $\nu$.

For a signed measure $\nu$, we use $\nu^+$ and $\nu^-$
to denote its positive and negative parts of $\nu$ respectively.
A singed measure $\nu$ is called smooth if both $\nu^+$ and $\nu^-$
are smooth.  For a signed smooth measure $\nu$, if
$A^{+}$ and $A^{-}$ are the PCAFs of $X^D$ with
Revuz measures $\nu^+$ and $\nu^-$ respectively,
the additive functional
$A:=A^+-A^{-}$ of is called the CAF of $X^D$ with
(signed) Revuz measure $\nu$.
When $\nu(dx) = c(x)dx$, $A_t$ is given by
$A_t=\int^t_0c(X^D_s)ds.$

We recall now the definition of the Kato class.

\medskip

\begin{defn}\label{df:3.1}
A signed smooth
measure $\nu$
is said to be in the class $\S_\infty (X^D)$
if for any $\eps>0$ there is a Borel
subset $K=K(\eps)$ of finite $| \nu|$-measure
and a constant $\delta = \delta (\eps) >0 $ such that
\begin{equation}\label{eqn:3.A}
\sup_{(x, z)\in (D\times D)\setminus d} \int_{D\setminus K}
\frac{G_D(x, y){G}_D(y, z) }{ {G}_D(x, z)} \, |\nu | (dy)\le\eps
\end{equation}
and for all measurable set $B\subset K$ with $| \nu |(B)<\delta$,
\begin{equation}\label{eqn:3.B}
\sup_{(x, z)\in (D\times D)\setminus d}
\int_B \frac{{G}_D(x, y){G}_D(y, z) }{{G}_D(x, z)} \, |\nu |
(dy) \le\eps.
\end{equation}
A function $q$ is said to be in the class
$\S_\infty (X^D)$,
if $q(x)dx$ is in $\S_\infty (X^D)$.
\end{defn}

\medskip

It follows from Proposition 7.1 of \cite{KS2}
and Theorem \ref{t:3G} above that $\K_{d,2}$ is
contained in $\S_\infty (X^D)$.
In fact, by Theorem \ref{t:3G2} we have the following result.
Recall that
$ g_1(x )=  G_D(x, z_0) \wedge C_1 $ and $
g_2(y )=  G_D(z_0, y) \wedge C_1$.

\medskip

\begin{prop}\label{p:3.2}
If a signed smooth measure $\nu$ satisfies
$$
\sup_{x\in D}\lim_{r\downarrow 0}
\int_{D\cap \{|x-y|\le r\}} \frac{g_1(y)}{g_1(x)} G_D(x,y)
|\nu | (dy) \, =\, 0
$$
and
$$
\sup_{x\in D}\lim_{r\downarrow 0}\int_{D\cap \{|x-y|\le r\}}
\frac{g_2(y)}{g_2(x)} G_D(y,x) |\nu | (dy)    \, =\, 0,
$$
then $\nu\in \S_\infty (X^D)$.
\end{prop}

\medskip

\pf This is a direct consequence of Theorem \ref{t:3G2}.
\qed

\medskip

In the remainder of this section, we will fix a signed measure
$\nu\in S_\infty(X^D)$ and we will use $A$ to denote the
CAF of $X^D$ with Revuz measure $\nu$. For simplicity, we will
use $e_A(t)$ to denote $\exp(A_t)$. The CAF $A$ gives rise to
a Schr\"odinger semigroup:
$$
Q^D_t f(x):= \E_x \left[ e_{A}(t)f(X^D_t) \right].
$$

The function $x\mapsto \E_x [e_{A}(\tau_D)]$ is called the gauge function
of $\nu$. We say $\nu$ is {\it gaugeable}
if $ \E_x \left[ e_{A} (\tau_D) \right]$ is finite
for some $x\in D$. In  the remainder of this
section we will assume that
$\nu$ is gaugeable.
It is shown in \cite{KS2}, by using the duality and the gauge theorems in
\cite{C2} and \cite{CS6}, that
the gauge function $x\mapsto \E_x [e_{A}(\tau_D)]$ is bounded on $D$
 (see section 7 in \cite{KS2}).

For $y\in D$, let $X^{D,y}$ denote the $h$-conditioned process obtained
from  $X^D$ with
$h(\cdot)= G_D(\cdot , y)$ and let $\E_x^y$ denote the expectation
for $X^{D,y}$ starting from $x\in D$. We will use $\tau^y_D$ to
denote the lifetime of $X^{D,y}$.
We know from  \cite{KS2} that $\E_x^y[ e_{A}(\tau_D^y)]$ is continuous in
$D \times D$ (also see Theorem 3.4 in \cite{CK2}) and
\begin{equation}\label{e:cg}
\sup_{(x, y)\in (D\times D)\setminus d} \E_x^y[
|A|_{\tau_D^y}]< \infty
\end{equation}
(also see \cite{C2} and \cite{CS6})
and therefore by Jensen's inequality
\begin{equation}\label{eqn:jen}
\inf_{(x, y)\in (D\times D)\setminus d} \E_x^y[ e_{A}(\tau_D^y)]>0,
\end{equation}
where $d$ is the diagonal of the set $D\times D$.
We also know from section 7 in \cite{KS2} that
\begin{equation}\label{e:green}
V_D(x, y):=\E_x^y[ e_{A}(\tau_D^y)]G_D(x, y)
\end{equation}
is the Green function of $\{Q^D_t\}$, that is, for any nonnegative
function $f$ on $D$,
$$
\int_D V_D(x, y) f(y) \,dy = \int_0^\infty
Q^D_t f (x) \, dt
$$
(also see Lemma 3.5 of \cite{C2}).
(\ref{e:cg})-(\ref{e:green}) and the continuity of
$\E_x^y[ e_{A}(\tau_D^y)]$
imply that
$V_D(x, y)$ is comparable to $G_D(x, y)$
and $V_D(x, y)$ is continuous on $(D\times D) \setminus d$.
Thus there exists a constant $c >0$
such that for every  $x, y, z \in D$,
\begin{equation}\label{3G_est4}
\frac{V_D(x,y) V_{D}(y,z)}
{ V_{D}(x,z)} \,\le\, c\, \frac{|x-z|^{d-2}}{|x-y|^{d-2} |y-z|^{d-2}}.
\end{equation}

\section{Two-sided heat kernel estimates for $\{Q^D_t\}$}

In this section, we will establish two-sided estimates for the heat kernel
of $Q^D_t$ in bounded $C^{1,1}$ domains.

Recall that a  bounded domain $D$ in $\R^d$ is said to be
a $C^{1,1}$
domain if there is a localization radius
$r_0>0$
 and a constant
$\Lambda >0$
such that for every $Q\in \partial D$, there is a
$C^{1,1}$-function $\phi=\phi_Q: \R^{d-1}\to \R$ satisfying $\phi (0)
= \nabla\phi (0)=0$,
$\| \nabla \phi  \|_\infty \leq \Lambda$,
$| \nabla \phi (x)-\nabla \phi (z)| \leq \Lambda
|x-z|$, and an orthonormal coordinate
system $y=(y_1, \cdots, y_{d-1}, y_d):=(\tilde y, y_d)$
such that
$ B(Q, r_0)\cap D=B(Q, r_0)\cap \{ y: y_d > \phi (\tilde y) \}$.

We will always assume in this section that $D$ is a bounded
$C^{1, 1}$ domain. Since
we will follow the method in \cite{KS} (see also \cite{Zq1}),
the proof of this section will be little sketchy.

First, we recall some results from \cite{KS}.
For every bounded $C^{1,1}$ domain $D$ and any $T>0$,
there exist positive constants $c_i$, $i=1, \dots, 4,$ such that
\begin{equation}\label{e:4.0}
C_{1}\psi_D(t, x, y)
t^{-\frac{d}2}
e^{-\frac{C_{2}|x-y|^2}t}
\le q^{D}(t, x, y)\le
C_{3}\psi_D(t, x, y)
t^{-\frac{d}2}e^{-\frac{C_{4}|x-y|^2}{t}}
\end{equation}
for all $(t, x, y)\in (0, T]\times D\times D$,
where
$$
\psi_D(t, x, y):=(1\wedge \frac{\rho_D(x)}{\sqrt t})
(1\wedge \frac{\rho_D(y)}{\sqrt t})$$
(see (4.27) in \cite{KS}).

For any $z\in \R^d$ and $ 0 < r \le 1 $, let
$$
D^z_r:=z+ rD, \quad
\psi_{D^z_r}(t, x, y):=(1\wedge \frac{\rho_{D^z_r}(x)}{\sqrt t})
(1\wedge \frac{\rho_{D^z_r}(y)}{\sqrt t}),
\quad (t, x, y)\in (0, \infty)\times D^z_r \times D^z_r
$$
where $\rho_{D^z_r}(x)$ is the distance between $x$ and
$\partial D^z_r $.
Then, for any $T>0$,
there exist positive constants  $t_0$ and $c_j, 5\le j\le 8,$
independent of $z$ and $r$
such that
\begin{equation}\label{e:4.1}
c_{5}t^{-\frac{d}{2}}\psi_{D^z_r}(t, x, y)
e^{-\frac{c_{6} |x-y|^2}{2t}}
\,\le \,q^{D^z_r}(t, x, y)\le c_{7}
t^{-\frac{d}{2}}\psi_{D^z_r}(t, x, y)e^{-\frac{c_{8}|x-y|^2}{2t}}
\end{equation}
for all $(t, x, y)\in (0, t_0\wedge (r^2T)]\times D^z_r\times D^z_r $
 (see (5.1) in \cite{KS}).
We will sometimes suppress the indices from $D^z_r$ when
there is no possibility of confusion.

For the remainder of this paper, we will assume that $\nu$ is in the
Kato class $\K_{d,2}$.
Using the estimates above and the joint continuity of the densities
$q^D(t,x,y)$ (Theorem 2.4 in \cite{KS1}), it is routine
(For example, see Theorem 3.17 \cite{CZ}, Theorem 3.1 \cite{BM} and
page 4669 in \cite{C2}.)
to show that $Q^D_t$ has a
jointly continuous density $r^D(t,\cdot,\cdot)$
(also see Theorem 2.4 in \cite{KS1}). So we have
\begin{equation}\label{e:dr}
\E_x\left[ e_A(t) f(X^D_t)   \right] = \int_D f(y) r^D(t,x,y) dy
\end{equation}
where $A$ is the CAF of $X^D$ with Revuz measure $\nu$ in $D$.

\medskip

\begin{thm}\label{t:dp}
The density $r^D(t,x,y)$ satisfies
the following equation
\begin{equation}\label{e:dp}
r^D(t,x,y)=q^D(t, x, y)+\int^t_0\int_{D}r^D(s,x,z)
 q^D(t-s, z, y)\nu(dz)ds
\end{equation}
for all $(t, x, y)\in (0, \infty)\times D\times D$.
\end{thm}
\pf
Recall that $A$ is the CAF of $X^D$ with Revuz measure $\nu$ in $D$ and
Let $\theta$ be the usual shift operator for Markov processes.

Since for any $t>0$
$$
e_A(t)=e^{A_t} = 1 + \int_{0}^{t} e^{A_t -A_s} dA_s=
1 + \int_{0}^{t} e^{A_{t-s} \circ \theta_s} dA_s,$$
We have
\begin{equation}\label{e:dp1}
\E_x\left[ e_A(t) f(X^D_t)   \right]
= \E_x\left[ f(X^D_t)   \right]
+ \E_x\left[ f(X^D_t) \int_{0}^{t}  e^{A_{t-s} \circ \theta_s} dA_s\right]
\end{equation}
for all $(t,x) \in (0, \infty) \times D$ and all bounded Borel-measurable
functions $f$ in D.

By the Markov Property and Fubini's theorem, we have
\begin{eqnarray*}
 \E_x\left[ f(X^D_t) \int_{0}^{t}   e^{A_{t-s} \circ \theta_s} dA_s\right]
&=&\int_{0}^{t}\E_x\left[ f(X^D_t)  e^{A_{t-s} \circ \theta_s} dA_s\right] \\
&=&\int_{0}^{t}\E_x\left[\E_{X^D_s}
\left[ f(X^D_{t-s})  e_A(t-s)\right]  dA_s\right].
\end{eqnarray*}
Thus by (\ref{eqn:revuz}) and (\ref{e:dr}),
\begin{equation}\label{e:dp2}
  \E_x\left[ f(X^D_t) \int_{0}^{t}   e^{A_{t-s} \circ \theta_s} dA_s\right]
=\int_D f(y)\int^t_0\int_{D}r^D(s,x,z)
 q^D(t-s, z, y)\nu(dz)ds dy.
\end{equation}
Since $r^D(s,\cdot, \cdot)$ and
 $q^D(t-s, \cdot, \cdot)$ are jointly continuous,
combining (\ref{e:dp1})-(\ref{e:dp2}), we have proved the theorem.
\qed

The proof of the next lemma is almost identical to that of
Lemma 3.1 in \cite{Zq2}. We omit the proof.

\medskip

\begin{lemma}\label{l:4.0}
For any $a>0$, there exists a positive constants $c$
depending only on $a$ and $d$ such that for
any $(t, x, y)\in (0, \infty)\times \R^d \times \R^d$,
\begin{eqnarray*}
&&\int^t_0\int_{\R^d} s^{-\frac{d}2}e^{-\frac{a|x-z|^2}{2s}}
(t-s)^{-\frac{d}2}
e^{-\frac{a|z-y|^2}{t-s}}|\nu|(dz)ds\\
&&\le ct^{-\frac{d}2}e^{-\frac{a|x-y|^2}{2t}}
\sup_{u\in \R^d}\int^t_0\int_{\R^d}s^{-\frac{d}2}e^{-\frac{a|u-z|^2}{4s}}
|\nu|(dz)ds
\end{eqnarray*}
and
\begin{eqnarray*}
&&\int^t_0\int_{\R^d} s^{-\frac{d+1}2}e^{-\frac{a|x-z|^2}{2s}}
(t-s)^{-\frac{d}2}
e^{-\frac{a|z-y|^2}{t-s}}|\nu|(dz)ds\\
&&\le ct^{-\frac{d+1}2}e^{-\frac{a|x-y|^2}{2t}}
\sup_{u\in \R^d}\int^t_0\int_{\R^d}s^{-\frac{d}2}e^{-\frac{a|u-z|^2}{4s}}
|\nu|(dz)ds
\end{eqnarray*}
\end{lemma}

\medskip

\begin{lemma}\label{l:4.1}
For any $a>0$, there exists a positive constant $c$
depending only on $a$ and $d$ such that for any $(t, x, y)
\in (0, \infty)\times D \times D$,
\begin{eqnarray}
&&\int^t_0\int_D(1\wedge \frac{\rho(x)}{\sqrt s})(1\wedge
\frac{\rho(z)}{\sqrt s})s^{-\frac{d}2}e^{-\frac{a|x-z|^2}{2s}}
(1\wedge\frac{\rho(y)}{\sqrt {t-s}})(t-s)^{-\frac{d}2}
e^{-\frac{a|z-y|^2}{t-s}}|\nu|(dz)ds\nonumber\\
&&\le c(1\wedge \frac{\rho(x)}{\sqrt t})(1\wedge
\frac{\rho(y)}{\sqrt t})t^{-\frac{d}2}e^{-\frac{a|x-y|^2}{2t}}
\sup_{u\in \R^d}\int^t_0\int_{\R^d}s^{-\frac{d}2}e^{-\frac{a|u-z|^2}{4s}}
|\nu|(dz)ds\label{e:4.4}
\end{eqnarray}
\end{lemma}
\pf
With Lemma \ref{l:4.0} in hand,
we can follow the proof of Theorem 2.1 (page 389-391)
in \cite{R1} to get the next lemma. So we skip the details.
\qed

\medskip

Recall that
$$
M^1_{\mu^i}(r)=\sup_{x\in \R^d}\int_{|x-y|\le r}\frac{|\mu^i|(dy)}
{|x-y|^{d- 1}} \quad
\mbox{and}
\quad
M^2_{\nu}(r)=\sup_{x\in \R^d}\int_{|x-y|\le r}\frac{|\nu|(dy)}
{|x-y|^{d- 2}}, \quad r >0, i=1 \cdots d.
$$

\medskip

\begin{thm}\label{t:4.2}
\begin{description}

\item{(1)}
For each $T>0$, there exist positive constants
$c_j, 1\le j\le 4,$
depending on $\mu$ and $\nu$ only via the rate
  at which $\max_{1 \le i \le d}
M^1_{\mu^i}(r)$ and $M^2_{\nu}(r)$ go  to zero
such that
\begin{equation}\label{est:5.0}
c_{1}t^{-\frac{d}{2}} \psi_D(t, x, y)
e^{-\frac{c_{2} |x-y|^2}{2t}}
\,\le \,r^{D}(t, x, y)\le c_{3}
t^{-\frac{d}{2}}\psi_D(t, x, y)e^{-\frac{c_{4}|x-y|^2}{2t}}
\end{equation}
\item{(2)}
There exist $T_1=T_1(D)>0$ such that for any $T>0$,
there exist positive constants  $t_1$ and $c_j, 5\le j\le 8,$
independent of $z$ and $r$ such that
\begin{equation}\label{est:5.1}
c_{5}t^{-\frac{d}{2}}\psi_{D^z_r}(t, x, y)
e^{-\frac{c_{6} |x-y|^2}{2t}}
\,\le \,r^{D^z_r}(t, x, y)\le c_{7}
t^{-\frac{d}{2}}\psi_{D^z_r}(t, x, y)e^{-\frac{c_{8}|x-y|^2}{2t}}
\end{equation}
for all $r \in (0,1]$ and $(t, x, y)\in (0, t_1\wedge (r^2(T \wedge T_1))  ]
\times D^z_r\times D^z_r $.
\end{description}
\end{thm}

\pf
We only give the proof of (\ref{est:5.1}). The proof of
(\ref{est:5.0}) is similar.
Fix $T>0$ and $z\in \R^d$.
Let $D_r :=D^z_r$, $\rho_{r}(x):=\rho_{D^z_r}(x)$ and
$\psi_{r}(t, x, y):=\psi_{D^z_r}(t, x, y)$.
We define $\tilde{I}_k(t,x,y)$ recursively for $k \ge 0$ and $(t, x , y)
\in (0, \infty) \times D \times D$:
\begin{eqnarray*}
I^r_0(t,x,y)&:=&q^{D_r}(t, x, y), \\
I^r_{k+1}(t,x,y)&:=&
\int^t_0\int_{D_r} I^r_k(s,x,z)  q(z) q^{D_r}(t-s, z, y)dzds.
\end{eqnarray*}
Then iterating the above gives
\begin{equation}\label{iter}
 r^{D_r}(t, x, y)= \sum_{k=0}^{\infty} I^r_k(t,x,y),
\quad  (t, x , y)\in (0, \infty) \times D_r \times D_r.
\end{equation}

Let
$$
N^2_{\nu}(t):=\sup_{u\in \R^d}\int^t_0\int_{\R^d}
       s^{-\frac{d}2}
e^{-\frac{|u-z|^2}{2s}}  |\nu|(dz) ds, \quad t >0 .
$$
It is well-known (See, for example, Proposition 2.1 in \cite{KS}.)
that for any $r >0$, there exist $c_1=c_1(d,r)$ and $c_2=c_2(d)$ such that
\begin{equation}\label{L_1}
N^2_{\nu}(t)\,\le\, (c_1t+c_2) M^2_{\nu}(r),
\quad \mbox{ for every } t \in (0,1).
\end{equation}

We claim that there exist positive constants $c_3$, $c_4$ and $A$
depending only on the
constants in (\ref{e:4.1}) and (\ref{e:4.4})
such that
for $k= 0, 1,\cdots$ and  $(t, x, y)\in (0, t_0\wedge (r^2T)]
\times D_r\times D_r$
\begin{equation}\label{inequal:J}
|I^r_k(t,x,y)|\, \le\, c_3 \,\psi_r(t, x, y)\,t^{-\frac{d}2}\,
e^{-\frac{|x-y|^2}{2t}}
 \left(c_4 N^2_{\nu}(\frac{2t}{A}) \right)^k, \quad 0< r \le 1.
\end{equation}
We will prove the above claim by induction.
By (\ref{e:4.1}), there exist constants $t_0$,
$c_3$ and $A$ such that
\begin{equation}\label{t:4.2_1}
|I^r_0(t,x,y)|=|q^{D_r}(t, x, y)| \,\le\, c_3\, \psi_r(t, x, y)\,
t^{-\frac{d}2}\,e^{-\frac{A|x-y|^2}{2t}}
\end{equation}
for $(t, x, y)\in (0, t_0\wedge (r^2T)]\times D_r\times D_r $.
On the other hand, by Lemma \ref{l:4.1}, there exists a
positive constant $c_5$ depending
only on $A$ and $d$ such that
\begin{eqnarray}
&&\int^t_0\int_{D_r}\psi_r(s, x, z)s^{-\frac{d}2}e^{-\frac{A|x-z|^2}{2s}}
(1\wedge\frac{\rho_r(y)}{\sqrt {t-s}})(t-s)^{-\frac{d}2}
e^{-\frac{A|z-y|^2}{t-s}} |\nu|(dz)  ds \nonumber\\
&&\le c_5\psi_r(t, x, y)t^{-\frac{d}2}e^{-\frac{A|x-y|^2}{2t}}
\sup_{u\in \R^d}\int^t_0\int_{\R^d}s^{-\frac{d}2}
e^{-\frac{A|u-z|^2}{4s}}
|\nu|(dz)ds. \label{t:4.2_2}
\end{eqnarray}
So there exists $c_6=c_6(d)$ such that
\begin{eqnarray*}
|I^r_1(t,x,y)| &\le& c_3^2 c_5\psi_r(t, x, y)t^{-\frac{d}2}
e^{-\frac{A|x-y|^2}{2t}}
\sup_{u\in \R^d}\int^t_0\int_{\R^d}s^{-\frac{d}2}
e^{-\frac{A|u-z|^2}{4s}}
|\nu|(dz)ds \\
&\le& c_3^2c_5 c_6 A^{\frac{d}{2}}
\psi_r(t, x, y)\,t^{-\frac{d}2}\,e^{-\frac{A|x-y|^2}{2t}}
N^2_{\nu}(\frac{2t}{A})
\end{eqnarray*}
for $(t, x, y)\in (0, t_0\wedge (r^2T)]\times D_r\times D_r $.
Therefore (\ref{inequal:J}) is true for $k=0,1$ with
$c_4 := c_3^2c_5 c_6 A^{\frac{d}{2}} $.
Now we assume (\ref{inequal:J}) is true up to $k$. Then by
(\ref{t:4.2_1})-(\ref{t:4.2_2}),
we have
\begin{eqnarray*}
&&|I^r_{k+1}(t,x,y)| \le
\int^t_0\int_{D_r} |I^r_k(s,x,z)| q^{D_r}(t-s, z, y)||\nu|(dz)ds\\
&& \le
\int^t_0\int_{D_r}  c_3\, \psi_r(s, x, z)\, s^{-\frac{d}2}
e^{-\frac{A|x-z|^2}{2s}}
 \left(c_4 \,  N^2_{\nu}(\frac{2s}{A}) \right)^k \\
&&\quad \times
  c_3 (1\wedge\frac{\rho_r(y)}{\sqrt {t-s}})(t-s)^{-\frac{d}2}
e^{-\frac{A|z-y|^2}{t-s}}|\nu|(dz)ds\\
&& \le c_3^2 \left(c_4 \, N^2_{\nu}(\frac{2t}{A}) \right)^k
\int^t_0\int_{D_r}  \psi_r(s, x, z)\,s^{-\frac{d}2}e^{-\frac{A|x-z|^2}{2s}}
 (1\wedge\frac{\rho_r(y)}{\sqrt {t-s}})\\
&&\quad \times
(t-s)^{-\frac{d}2}
e^{-\frac{M|z-y|^2}{t-s}}|\nu|(dz)ds \\
&& \le c_3^2 \left(c_4 \, N^2_{\nu}(\frac{2t}{A}) \right)^k
c_5 c_6 A^{\frac{d}{2}}
\psi_r(t, x, y)t^{-\frac{d}2}e^{-\frac{A|x-y|^2}{2t}} N^2_{\nu}(\frac{2t}{A})
\\
&& \le
c_3 \psi_r(t, x, y)t^{-\frac{d}2}e^{-\frac{A|x-y|^2}{2t}}
\left(c_4 \, N^2_{\nu}(\frac{2t}{A}) \right)^{k+1}.
\end{eqnarray*}
So the claim is proved.

Choose $t_1<(1 \wedge t_0)$ small so that
\begin{equation}\label{N_1}
c_4 \, N^2_{\nu}(\frac{2t_1}{A})  <\frac12.
\end{equation}
By (\ref{L_1}), $t_1$ depends on $\nu$ only via the
rate at which $M^2_{\nu}(r)$ goes to zero.
(\ref{iter}) and (\ref{inequal:J}) imply that
for $  (t, x , y)\in (0,t_1 \wedge (r^2T)] \times D_r \times D_r$
\begin{equation}\label{qup}
 r^{D_r}(t, x, y)\,\le\, \sum_{k=0}^{\infty} |I^r_k(t,x,y)|
\,\le\,   2 c_3 \psi_r(t, x, y)t^{-\frac{d}2}e^{-\frac{A|x-y|^2}{2t}}.
\end{equation}

Now we are going to prove the lower estimate of
$r^{D_r}(t, x, y)$.
Combining  (\ref{iter}),  (\ref{inequal:J}) and (\ref{N_1}) we have
for every $  (t, x , y)\in (0,t_1 \wedge (r^2T)] \times D_r \times D_r$,
\begin{eqnarray*}
&& |r^{D_r}(t, x, y)- q^{D_r}(t, x, y)|\le
\sum_{k=1}^{\infty} |I^r_k(t,x,y)|
\,\le\, c_3 c_4 \,  N^2_{\nu}(\frac{2 t_1 }{A})
\psi_r(t, x, y)t^{-\frac{d}2}e^{-\frac{A|x-y|^2}{2t}}.
\end{eqnarray*}
Since there exist $c_7$ and $c_8 \le 1$ depending
on $T$ such that
$$
q^{D_r}(t, x, y) \ge 2c_8 \psi_r(t, x, y)t^{-\frac{d}2}
e^{-\frac{c_7|x-y|^2}{2t}},
$$
we have for $|x-y| \le   {\sqrt t}$ and $(t, x, y)\in
(0,   t_1 \wedge (r^2T) ]\times D\times D$,
\begin{equation}\label{lower3}
r^{D_r}(t, x, y)
\,\ge\, \left(2c_8e^{-2 c_7} -  c_3 c_4 \, N^2_{\nu}(\frac{2 t_1 }{A})
\right) \psi(t, x, y)t^{-\frac{d}2}.
\end{equation}
Now we choose $t_2 \le t_1$ small so that
\begin{equation}\label{N_3}
c_3 c_4 \,  N^2_{\nu}(\frac{2 t_2 }{A}) \,<\,
c_8e^{-2 c_7}.
\end{equation}
Note that $t_2$ depends on $\nu$ only via the
rate at which $M^2_{\nu}(r)$ goes to zero.
So for   $(t, x, y)\in
(0,   t_2 \wedge (r^2T) ]\times D\times D$ and
$|x-y| \le  {\sqrt t}$, we have
\begin{equation}\label{lower2}
 r^{D_r}(t, x, y)
\ge c_8e^{-2c_7} \psi_r(t, x, y)t^{-\frac{d}2}.
\end{equation}

It is easy to check (see pages 420--421
of \cite{Zq3}) that there exists a positive constant
$T_0$ depending only on the characteristics of
the bounded $C^{1, 1}$ domain $D$ such that
for any $\hat t\le T_0$ and $x, y\in D$ with $\rho_D(x)\ge \sqrt{\hat t},
\rho_D(y)\ge \sqrt{\hat t}$, one can find a arclength-parameterized
curve $l\subset D$ connecting $x$ and $y$ such that the length
$|l|$ of $l$ is equal to $\lambda_1|x-y|$ with $\lambda_1\le
\lambda_0$, a constant depending only on the characteristics of
the bounded $C^{1, 1}$ domain $D$. Moreover, $l$ can be chosen so that
$$
\rho_D(l(s))\ge \lambda_2\sqrt{\hat t}, \quad s\in [0, |l|]
$$
for some positive constant $\lambda_2$ depending only on the
characteristics of the bounded $C^{1, 1}$ domain $D$.
Thus for
any $t=r^2\hat t\le r^2T_0$ and $x, y\in D_r$ with $\rho_r(x)\ge \sqrt{ t},
\rho_r(y)\ge \sqrt{ t}$, one can find a arclength-parameterized
curve $l\subset D_r$ connecting $x$ and $y$ such that the length
$|l|$ of $l$ is equal to $\lambda_1|x-y|$ and
$$
\rho_r(l(s))\ge \lambda_2\sqrt{ t}, \quad s\in [0, |l|].
$$

Using this fact and (\ref{lower2}), and following the proof of
Theorem 2.7 in \cite{FS}, we can show that there exists a
positive constant $c_9$ depending only on $d$ and the
characteristics of the bounded $C^{1, 1}$ domain $D$
such that
\begin{equation}\label{lower4}
r^{D_r}(t, x, y)
\ge \frac12c_8e^{-2c_7} \psi_r(t, x, y)t^{-\frac{d}2}e^{-\frac{c_9|x-y|^2}t}
\end{equation}
for all $t\in (0, t_2 \wedge r^2(T \wedge T_0)]$ and
$x, y\in D_r$ with $\rho_r(x)\ge \sqrt{t},
\rho_r(y)\ge \sqrt{t}$.

It is easy to check that there exists a positive constant
$T_1 \le T_0 $ depending only on the characteristics of
the bounded $C^{1, 1}$ domain $D$ such that
for $\hat t \le T_1 $ and arbitrary $x, y\in D$,
one can find $x_1, y_1\in D$ be such that
$\rho_D(x_1)\ge \sqrt{\hat t}$, $\rho_D(y_1)\ge \sqrt{\hat t}$
and $|x-x_0|\le \sqrt{\hat t}$, $|y-y_0|\le \sqrt{\hat t}$.
Thus for
any $t=r^2\hat t\le r^2T_1$
 and arbitrary $x, y\in D_r$,
one can find $x_1, y_1\in D_r$ be such that
$\rho_r(x_1)\ge \sqrt{t}, \rho_r(y_1)\ge \sqrt{t}$
and $|x-x_0|\le \sqrt{t}, |y-y_0|\le \sqrt{t}$. Now
 Using (\ref{lower3}) and (\ref{lower4})
one can repeat the last paragraph of the proof of
Theorem 2.1 in \cite{R1} to show that there exists
a positive constant $c_{10}$ depending only on $d$ and
the characteristics of the bounded $C^{1, 1}$ domain $D$
such that
\begin{equation}\label{lower5}
r^{D_r}(t, x, y)
\ge c_8c_{10}e^{-2c_7} \psi_r(t, x, y)t^{-\frac{d}2}e^{-\frac{2c_9|x-y|^2}t}
\end{equation}
for all $(t, x, y)\in (0, t_2 \wedge r^2(T \wedge T_1)]\times D_r\times D_r$.

Using (\ref{e:4.0}) instead of  (\ref{e:4.1})
The proof of (\ref{est:5.0}) up to $t \le t_3$ for some $t_3$
depending on $T$ and $D$
is similar (and simpler) to the proof of (\ref{est:5.1}).
To prove (\ref{est:5.0}) for a general $T>0$, we can
apply the Chapman-Kolmogorov equation and use the argument
in the proof of Theorem 3.9 in \cite{So}. We omit the details.
\qed

\medskip

\begin{remark}\label{r1}
{\rm Theorem \ref{t:4.2} (2) will be used in \cite{KS5} to prove
parabolic Harnack inequality, parabolic
boundary Harnack inequality and  the intrinsic ultracontractivity for the
semigroup $Q_t^D$. }
\end{remark}

\section{Uniform 3G type estimates for small Lipschitz domains}

Recall that $r_1>0$ is the constant from (\ref{e:Green_B}) and $r_3>0$
is the constant from Theorem \ref{HP}.
The next lemma is a scale invariant version of Lemma \ref{G:HP}.
The proof is similar to the proof
of Lemma \ref{G:HP}.

\medskip

 \begin{lemma}\label{G:HP1}
There exists $c=c(d, \mu) >0$ such that for  every $r \in
(0 , r_1 \wedge r_3]$, $Q \in \R^d$ and
open subset $U$ with $B(z, l) \subset U \subset B(Q,r)$,
we have for every $ x \in
U \setminus  \overline{B(z, l)}$
\begin{equation}\label{e:GH11}
\sup_{y \in B(z, l/2)} G_{U}(y, x)
\le c \inf _{y \in B(z, l/2)} G_{U}(y, x)
\end{equation}
and
\begin{equation}\label{e:GH21}
\sup_{y \in B(z, l/2)} G_{U}(x, y)
\le c \inf _{y \in B(z, l/2)} G_{U}(x, y)
\end{equation}
\end{lemma}

\pf (\ref{e:GH11}) follows
from Theorem  \ref{HP}. So we only need to show (\ref{e:GH21}).
Since $r < r_1$, by (\ref{e:Green_B}),
there exists $c=c(d) >1$ such that
 for every $x, w \in \overline{B(z, \frac{3l}4)}$
$$
c^{-1}\, \frac{1}{|w-x|^{d-2}} \,\le\,   G_{B(z,l)}  (w,x)
\,\le\,  G_{U}  (w,x) \,\le\, G_{B(Q, r )}  (w,x) \,\le\, c \,
\frac{1}{|w-x|^{d-2}}.
$$
Thus for $ w \in \partial {B(z, \frac{3l}4)}$ and $ y_1, y_2 \in
{B(z, \frac{l}2)}$, we have
\begin{equation}\label{e:pp31}
G_{U}  (w,y_1) \,\le\, c \,\left(\frac{|w-y_2|}{|w-y_1|}\right)^{d-2}
\frac{1}{|w-y_2|^{d-2}} \,\le\, 4^{d-2}\, c^2\, G_{U}  (w,y_2).
\end{equation}
On the other hand, from (\ref{e:GH}), we have
\begin{equation}\label{e:pp41}
 G_{U}(x,y) \, = \, \E_x\left[G_U(X_{T_{B(z,\frac{l}2)}},y)\right],
\quad y \in B(z, \frac{l}{2})
\end{equation}
Since $X_{T_{B(z,\frac{3l}4)}} \in \partial {B(z, \frac{3l}4)}$,
combining (\ref{e:pp31})-(\ref{e:pp41}), we get
$$
G_{U}(x,y_1)  \,\le\,  4^{d-2}\, c^2\,
\E_x\left[G_{U}(X_{T_{B(z,\frac{3l}4)}},y_2)\right]
\,=\, 4^{d-2}\, c^2 \, G_{U}(x,y_2)  , \quad y_1, y_2 \in B(z,\frac{l}2)
$$
\qed

\medskip

In the remainder of this section, we fix a bounded
Lipschitz domain $D$ with
characteristics $(R_0, \Lambda_0)$.
For every $Q\in \partial D$ we put
$$
\Delta_Q(r) :=\{ y \mbox{ in } CS_Q:
\phi_Q (\tilde y)+r >y_d > \phi_Q (\tilde y),\,
 |\tilde y| < r \}
$$
where $CS_Q$ is the coordinate
system with origin at $Q$ in the definition
of Lipschitz domains and $\phi_Q$ is the Lipschitz
function there. Define
\begin{equation}\label{e:r_4}
r_5:=\frac{R_0}{\sqrt{ 1 +\Lambda_0^2}+1}\wedge r_1 \wedge r_3.
\end{equation}
If $z \in \overline{\Delta_Q(r)}$ with $r\le r_5$,
we have
$$
|Q-z| \,\le\,
|( \tilde z, \phi_Q (\tilde z))-(\tilde{Q}, 0)| +r
\,\le\, (\sqrt{ 1 +\Lambda_0^2}+1)r   \le R_0.
$$
So  $\overline{\Delta_Q(r)} \subset B(Q, R_0)\cap D$.

For any Lipschitz function $\psi: \R^{d-1} \to \R$ with Lipschitz constant
$\Lambda_0$, let
$$
\Delta^\psi  \,:=\,\left\{ y: r_5 >
y_d - \psi (\tilde y) >0,\,
 |\tilde y | < r_5  \right\}.
$$
so that $\Delta^\psi \,\subset\,  B(0, R_0).$
We observe that, for any Lipschitz function $\varphi: \R^{d-1} \to \R$
with the Lipschitz constant $\Lambda$,
its dilation $\varphi_r (x) := r \varphi (x/r)$ is also Lipschitz
with the same Lipschitz constant $\Lambda_0$.
For any $r>0$, put
$\eta=\frac{r}{r_5}$
and $\psi=(\phi_Q)_{\eta}$. Then it is easy to see that
for any $Q \in \partial D$ and $r \le r_5$,
$$
\Delta_Q(r)=\eta \Delta^\psi.
$$
Thus by choosing appropriate constants $\Lambda_1>1$, $R_1 < 1$
and $d_1 >0$, we can say  that
for every $Q \in \partial D$ and  $r \le r_5$,
the $\Delta_Q(r)$'s are bounded Lipschitz domains with
the  characteristics $(rR_1, \Lambda_1)$ and the diameters of
$\Delta_Q(r)$'s are less than $rd_1$.
Since $r_5 \le r_1 \wedge r_3$, Lemma \ref{G:HP1} works for
$G_{\Delta_Q(r)}(x,y)$ with
$Q \in \partial D$ and  $r \le r_5$. Moreover, we can
restate the scale invariant boundary Harnack principle in the
following way.

\medskip

\begin{thm}\label{BHPi}
There exist  constants $M_3, c >1$  and $s_1 >0 $, depending on $\mu$,
$\nu$  and $D$
such that for every $Q \in \partial D$, $r < r_5$,
$s< r s_1$, $w \in \partial \Delta_Q(r)$ and any
nonnegative functions
$u$ and $v$ which are harmonic with respect to $X^D$ in
$\Delta_Q(r) \cap B(w, M_3 s)$
and vanish continuously on
$\partial \Delta_Q(r) \cap B(w, M_3s)$, we have
\begin{equation}\label{e:BHP}
\frac{ u(x)}{ v(x)}
\, \le \,c\, \frac{ u(y)}{ v(y)}  \quad  \mbox{ for any }
x,y \in \Delta_Q(r) \cap B(w, s).
\end{equation}
\end{thm}

\medskip

In the remainder of this section we will fix the above constants
$r_5$, $M_3$, $s_1$, $\Lambda_1$, $R_1$ and $d_1 >0$,
and consider the Green functions of $X$ in
$\Delta_Q(r)$ with $Q \in \partial D$ and $r>0$.
We will prove a scale invariant 3G type estimates for these Green
functions for small $r$.
The main difficulties of the scale invariant 3G type estimates
for $X$ are the facts that $X$ does not have rescaling property
and that the Green function
$G_{\Delta_Q(r)}(x, \, \cdot\,)$ is not
harmonic for $X$. To overcome these difficulties,
we first establish some results for the Green functions of $X$ in
$\Delta_Q(r)$ with $Q \in \partial D$ and  $r$ small.

Let $\delta^Q_r (x):=$ dist$(x, \partial \Delta_Q(r))$. Using
Lemma \ref{G:HP1} and a Harnack
chain argument, the proof of the next lemma is almost identical
to the proof of Lemma 6.7
 in \cite{CZ}. So we omit the proof.

\medskip

\begin{lemma}\label{l:Green_L1}
For any given $c_1>0$,
there exists $c_2=c_2(D, c_1, \mu)>0$ such that for every
$Q \in \partial D$, $r< r_5$, $|x-y| \le  c_1 (\delta^Q_r (x)
\wedge \delta^Q_r (y))$,
we have
$$
G_{\Delta_Q(r)}(x,y) \,\ge\, c_2 \,|x-y|^{-d+2}.
$$
\end{lemma}

\medskip

Recall that $M_3>0$ and $s_1 >0 $ are the constants
from Theorem \ref{BHPi}.  Let $M_4 :=
2(1+M_3)\sqrt{ 1 +\Lambda_1^2}+2$ and
$
R_4:=R_1/M_4.
$
The next lemma is a scale invariant version of Lemma \ref{l:BHPi}.
The proof is similar to the proof
of Lemma \ref{l:BHPi}. We spell out the details for the reader's convenience.

\medskip

\begin{lemma}\label{l:BHPi1}
There exists  constant $ c >1$
such that for every $Q \in \partial D$, $r < r_5$, $s < rR_4$,
$w \in \partial \Delta_Q(r)$ and any
nonnegative functions
$u$ and $v$ which are harmonic in $\Delta_Q(r) \setminus B(w,s)$
and vanish continuously on
$\partial \Delta_Q(r) \setminus B(w,s)$, we have
\begin{equation}\label{e:l:BHP1}
\frac{ u(x)}{ u(y)}
\, \le \,c\, \frac{ v(x)}{ v(y)}  \quad  \mbox{ for any }
x, y \in \Delta_Q(r) \setminus B(w, M_4s).
\end{equation}
\end{lemma}

\pf
We fix a point $Q$ on $\partial D$, $r < r_5$, $s < rR_4$ and $w \in \partial
\Delta_Q(r)$ throughout this proof.
Let
\begin{eqnarray*}
\Delta^s &:=&\left\{ y\mbox{ in } CS_w
:\, \varphi_w (\tilde y)+ 2s
\, > \,y_d\, >\, \varphi_w (\tilde y),\,\,
|\tilde y| < 2(M_3+1)s \right\},\\
\partial_1 \Delta^s   &:=& \left\{ y\mbox{ in } CS_w
:\, \varphi_w (\tilde y)+ 2s
 \,\ge\, y_d \,>\, \varphi_w
(\tilde y), \,\,|\tilde y| = 2(M_3+1)s \right\},\\
\quad \partial_2 \Delta^s
&:=& \left\{ y\mbox{ in } CS_w
: \,\varphi_w (\tilde y)+  2s
\,   =\, y_d
, \,\,|\tilde y| \le  2(M_3+1)s  \right\},
\end{eqnarray*}
where $CS_w$ is the coordinate
system with origin at $w$ in the definition
of the Lipschitz domain $\Delta_Q(r)$ and $\varphi_w$ is the Lipschitz
function there.
If $z \in \overline{\Delta^s}$,
$$
|w-z| \,\le\,
|( \tilde z, \varphi_w (\tilde z))-(\tilde{z}, 0)| +2s
\,\le\, 2s(1+M_3)\sqrt{ 1 +\Lambda^2}+2s = M_4s\,\le\, rR_1.
$$
So  $\overline{\Delta^s} \subset B(Q, M_4s)\cap D \subset  B(Q, rR_1)\cap D $.
For $|\tilde y| = 2(M_3+1)s$, we have $|(\tilde y,
\varphi_w (\tilde y))| > s$. So
$u$ and $v$ are harmonic with respect to $X$ in $ \Delta_Q(r)  \cap
B((\tilde y,  \varphi_w (\tilde y)), 2M_3s)$
and vanish continuously on
$\partial  \Delta_Q(r) \cap B((\tilde y,  \varphi_w (\tilde y)),2M_3s)$ where
$|\tilde y| =2(M_3+1)s$.
Therefore by Theorem \ref{BHPi},
\begin{equation}\label{B_111}
\frac{ u(x)}{ u(y)}
\, \le \,c\, \frac{ v(x)}{ v(y)}  \quad  \mbox{ for any }
x, y \in    \partial_1 \Delta^s \mbox{ with } \tilde x = \tilde y.
\end{equation}

Since dist$(  \Delta_Q(r) \cap B(w,s),\partial_2 \Delta^s) > c s$
for some $c_1=c_1(D)$,
if $x \in  \partial_2 \Delta^s$, the Harnack inequality
(Theorem \ref{HP}) and a
Harnack chain argument give that
there exists constant $c_2 >1$ such that
\begin{equation}\label{B_121}
c_2^{-1}\,<\,\frac{u(x)}{u(y)},\,\,\frac{v(x)}{v(y)} \,<\,c_2.
\end{equation}
In particular, (\ref{B_121}) is true with $x=
x_s:=( \tilde x, \varphi_w(\tilde x)+2s)$, which is also in
$\partial_1 \Delta^s$. Thus (\ref{B_111}) and (\ref{B_121})
imply that
\begin{equation}\label{B_141}
c_3^{-1} \frac{ u(x)}{ u(y)}
\, \le \, \frac{ v(x)}{ v(y)}  \le \,c_3\,\frac{ u(x)}{ u(y)},
\quad x,y \in  \partial_1 \Delta^s
 \cup  \partial_2 \Delta^s
\end{equation}
for some $c_3>1$.
Now, by applying the maximum principle (Lemma 7.2 in \cite{KS})
twice ($x$ and $y$),
(\ref{B_141}) is true for every $x \in  \Delta_Q(r)
\setminus  \Delta^s$.
\qed

\medskip

Combining Theorem \ref{BHPi} and Lemma \ref{l:BHPi1}, we get
the following as a corollary.

\medskip

\begin{corollary}\label{C:Green_L}
There exists  constant $ c >1$
such that for every $Q \in \partial D$, $r < r_5$,
$w \in \partial (\Delta_Q(r)) $,    and    $s < rR_4$,  we have for
$x, y \in \Delta_Q(r) \setminus B(w, M_4s)$ and $z_1, z_2
\in \Delta_Q(r) \cap B(w, s)$
\begin{equation}\label{e:CG_11}
\frac{ G_{\Delta_Q(r)}(x,z_1)}{ G_{\Delta_Q(r)} (y,z_1)}
\, \le \,c\, \frac{ G_{\Delta_Q(r)} (x,z_2) }{ G_{\Delta_Q(r)} (y,z_2)}
\quad
\mbox{and}
\quad
\frac{ G_{\Delta_Q(r)}(z_1,x)}{ G_{\Delta_Q(r)} (z_1,y)}
\, \le \,c\, \frac{ G_{\Delta_Q(r)} (z_2,x) }{ G_{\Delta_Q(r)} (z_2,y)}.
\end{equation}
\end{corollary}

\medskip

\begin{corollary}\label{C:c_L1}
For any given $N\in(0, 1)$, there exists  constant $ c=c(N, M_4, D) >1$
such that for every $Q \in \partial D$, $r < r_5$,
$w \in \partial (\Delta_Q(r)) $ and  $s < rR_4$,  we have
\begin{equation}\label{e:CG_31}
G_{\Delta_Q(r)}(x,z_1)
\, \le \,c\,  G_{\Delta_Q(r)} (x,z_2)
\quad
\mbox{and}\quad
 G_{\Delta_Q(r)}(z_1,x)
\, \le \,c\,  G_{\Delta_Q(r)} (z_2,x)
\end{equation}
for $x \in \Delta_Q(r) \setminus B(w, M_4s)$ and $z_1, z_2
\in \Delta_Q(r) \cap B(w, s)$
with $B(z_2, Ns) \subset \Delta_Q(r) \cap B(w, s)$.
\end{corollary}

\pf  Fix $Q \in \partial D$, $r < r_5$,
$w \in \partial (\Delta_Q(r)) $ and  $s < rR_4$.
Recall from the proof of Lemma \ref{l:BHPi1} that $CS_w$ is the coordinate
system with origin at $w$ in the definition
of the Lipschitz domain $\Delta_Q(r)$. Let $\overline{y}
:=(\tilde 0, M_4 s)$. By
(\ref{G_bd}),
$$
G_{\Delta_Q(r)} (\overline{y},z_1) \,\le\, c_1\,
|y-z_2|^{-d+2} \,\le\, c_2 s^{-d+2}
\quad
\mbox{and}
\quad
 G_{\Delta_Q(r)} (z_1,\overline{y}) \,\le\, c_1\,
|y-z_2|^{-d+2} \,\le\, c_2\,s^{-d+2},
$$
for some constants $c_1, c_2>0$.

Note that, since $\Delta_Q(r)$'s are bounded Lipschitz domains with
the  characteristics $(rR_1, \Lambda_1)$ and $s < rR_4 $,
it is easy to see that
there exists a positive constant $c_3$ such that
$\rho_r^Q(\overline{y}) \ge c_3 M_4s$ and
$\rho_r^Q(z_2) \ge N s$. Thus by Lemma \ref{l:Green_L1},
$$ G_{\Delta_Q(r)} (y,z_2)
\,\ge\, c_4\,
|y-z_2|^{-d+2}
\,\ge \,c_5\, s^{-d+2}
\quad
\mbox{and}
\quad
 G_{\Delta_Q(r)} (z_2,y) \,\ge\, c_4\,
|y-z_2|^{-d+2}
\,\ge \,c_5\, s^{-d+2}
$$
for some constants $c_4, c_5>0$.

Now apply (\ref{e:CG_11}) with $y=\overline{y}$ and get
$$
G_{\Delta_Q(r)}(x,z_1) \,\le\, c_6\,G_{\Delta_Q(r)} (x,z_2)
\quad
\mbox{and}\quad
G_{\Delta_Q(r)}(z_1,x) \,\le\, c_6\, G_{\Delta_Q(r)} (z_2,x),
$$
for some $c_6>1$.
\qed

\medskip

With lemma \ref{G:HP1}, Corollary \ref{C:Green_L} and
Corollary \ref{C:c_L1} in hand, one can follow either the argument
 in Section 2 of this paper or the argument on page 170-173
of \cite{CZ}. So we skip the details.

\medskip

\begin{thm}\label{t:3G1}
There exists a constant $c >0$
 such that for every $Q \in \partial D$, $r < r_5$ and  $x, y, z
\in \Delta_Q(r)$,
\begin{equation}\label{3G_est1}
 \frac{G_{\Delta_Q(r)}(x,y) G_{\Delta_Q(r)}(y,z)}
{ G_{\Delta_Q(r)}(x,z)} \,\le\, c\,\left(|x-y|^{-d+2} +
|y-z|^{-d+2}\right).
\end{equation}
\end{thm}

\section{Boundary Harnack principle for the Schr\"{o}dinger operator
of $X^D$ in bounded Lipschitz domains}

Recall that  $\nu$
belongs to the Kato class $\K_{d,2}$ and
$A$ is continuous additive
functional associated with $\nu|_D$. We also
recall $
e_{A}(t)=\exp( A_t)
$ and the Schr\"odinger semigroup
$$ Q^D_t f(x)= \E_x \left[ e_{A}(t)f(X^D_t) \right].
$$

Using
the Martin representation for Schr\"{o}dinger operators
(Theorem 7.5 in \cite{CK2}) and
the uniform 3G estimates (Theorem \ref{t:3G1}),
we will prove the boundary Harnack principle for
the Schr\"{o}dinger operator
 of diffusions with measure-valued drifts in bounded Lipschitz domains.
In the remainder of this section, we fix a bounded
Lipschitz domain $D$ with its
characteristics $(R_0, \Lambda_0)$. Recall
$$
\Delta_Q(r) =\{ y \mbox{ in } CS_Q:
\phi_Q (\tilde y)+r >y_d > \phi_Q (\tilde y),\,
|\tilde y| < r \},
$$
where $CS_Q$ is the coordinate
system with origin at $Q \in \partial D$ in the definition of
Lipschitz domains and $\phi_Q$ is the Lipschitz function there.
We also recall that $r_5$ is the constant from (\ref{e:r_4}) and
that the diameters of
$\Delta_Q(r)$'s are less than $rd_1$.

For $Q \in \partial D$, $r < r_5$ and $y \in \Delta_Q(r)$, let
$X^{Q,r,y}$ denote the $h$-conditioned process obtained
from  $X^{\Delta_Q(r)} $ with
$h(\cdot)= G_{\Delta_Q(r)}(\cdot , y)$ and let $\E_x^{Q,r,y}$
denote the expectation
for $X^{Q,r,y}$ starting from $x\in  \Delta_Q(r)$.
Now define the conditional gauge function
$$
u^Q_r(x,y):= \E_x^{Q,r,y}\left[ e_{A^\nu}
\left(\tau_{\Delta_Q(r)}^y\right)\right].
$$
By Theorem \ref{t:3G1},
\begin{eqnarray*}
\E_x^{Q,r,y}\left[A\left(\tau_{\Delta_Q(r)}^y\right)\right]
&\le&
\int_{\Delta_Q(r)} \frac{G_{\Delta_Q(r)}(x,a) G_{\Delta_Q(r)}(a,y)}
{ G_{\Delta_Q(r)}(x,y)} \nu(da) \\
&\le& c\,\int_{\Delta_Q(r)} \left(|x-a|^{-d+2} +
|a-y|^{-d+2}\right) \nu(da),\quad
r <r_5.
\end{eqnarray*}
Since the above constant is independent of $r < r_5$, we have
$$
\sup_{x,y \in \Delta_Q(r)}  \E_x^{Q,r,y}\left[A
\left(\tau_{\Delta_Q(r)}^y\right)\right]
\,\le \,
 c\,\sup_{x\in \R^d}\int_{|x-a|\le rd_1}\frac{|\nu|(da)}
{|x-a|^{d-2 }} \,=\,cM^2_{\nu}(rd_1) < \infty ,\quad
r <r_5, Q \in \partial D.
$$
Thus $\nu \in \S_\infty(X^{\Delta_Q(r)})$ for every $r < r_5$
and there exists  $r_6 \le
r_5$ such that $$
 \sup_{x,y \in\Delta_Q(r)} \E_x^{Q,r,y}\left[A\left(
\tau_{\Delta_Q(r)}^y\right)\right]
\,\le \, \frac12 ,\quad
r <r_6, Q \in \partial D.
$$
Hence by Khasminskii's lemma,
$$
\sup_{x,z \in\Delta_Q(r)}  u^Q_r(x,y) \,\le\, 2, \qquad  r <r_6,\, Q \in \partial D.
$$
By Jensen's inequality, we also have
$$
\inf_{x,z \in\Delta_Q(r)}  u^Q_r(x,y) \,>\, 0, \qquad  r <r_6,\, Q \in \partial D.
$$
Therefore, we have proved the following lemma.

\medskip

\begin{lemma}\label{l:k_b}
For $r <r_6$, $\nu |_{\Delta_Q(r)} \in \S_\infty(X^{\Delta_Q(r)})$
and $\nu |_{\Delta_Q(r)}$ is gaugeable.
Moreover, there exists a constant $c$ such that $ c^{-1} \,\le\,
u^Q_r(x,y)  \,\le\,c$
for $x,y \in \Delta_Q(r)$ and $r <r_6$.
\end{lemma}

\medskip

\begin{thm}\label{BHP4} (Boundary Harnack principle)
Suppose $D$ be a bounded Lipschitz domain in $\R^d$ with the
Lipschitz characteristic $(R_0, \Lambda_0)$ and let
$M_5:=(\sqrt{ 1 +\Lambda_0^2}+1)$.
Then there exists $N>1$ such that for any $r\in (0, r_6)$ and
$Q\in \partial D$,
there exists a constant $c>1$
such that for any nonnegative functions $u, v$
which are $\nu$-harmonic in $D \cap B(Q, rM_5)$ with respect to $X^D$
and  vanish continuously on
$\partial D \cap B(Q, rM_5)$, we have
$$
\frac{ u(x)}{ v(x)}
\, \le \,c\, \frac{ u(y)}{ v(y)}  \quad  \mbox{ for any }
x,y \in D\cap B(Q, \frac{r}{N}).
$$
\end{thm}

\pf
Note that, with $M_5=(\sqrt{ 1 +\Lambda_0^2}+1)$,
$\Delta_Q(r) \subset D \cap B(Q, M_5r)$. So
$u, v$ are $\nu$-harmonic in $\Delta_Q(r)$.
For the remainder of the proof, we fix  $Q\in \partial D$, $r\in (0, r_5)$
and a point $x^Q_r \in \Delta_Q(r)$.  Let
$$
M(x, z)
:= \lim_{U \ni y  \to z} \frac{G_{U}(x,y)}
{G_{U}(x^Q_r,y)}, \quad  K(x, z)
:= \lim_{U \ni y  \to z} \frac{V_{U}(x,y)}
{V_{U}(x^Q_r,y)}.
$$
Since $u, v$ are $\nu$-harmonic with respect to $ X^{\Delta_Q(r)}$,
by Theorem 7.7 in \cite{CK2} and our Lemma \ref{l:k_b}, there
exist finite measures $\mu_1$
and $\nu_1$ on $\partial U$
such that
$$
u(x)= \int_{\partial \Delta_Q(r)} K(x, z) \mu_1(dz) \quad \mbox{and}\quad
v(x)=\int_{\partial \Delta_Q(r)} K(x, z) \nu_1(dz), \quad x \in \Delta_Q(r).
$$
Let $$
u_1(x):= \int_{\partial \Delta_Q(r)} M(x, z) \mu_1(dz) \quad \mbox{and}\quad
v_1(x):=\int_{\partial \Delta_Q(r)} M(x, z) \nu_1(dz), \quad x \in \Delta_Q(r).
$$
By Theorem 7.3 (2) in \cite{CK2} and our Lemma \ref{l:k_b}, we have for every $x \in U$
\begin{eqnarray*}
&&\frac{u(x)}{v(x)}\,=\,\frac{\int_{\partial\Delta_Q(r)} K(x, z) \mu_1(dz)}
{\int_{\partial \Delta_Q(r)} K(x, z) \nu_1(dz)}
\,\le\,c_1^2\frac{\int_{\partial \Delta_Q(r)} M(x, z) \mu_1(dz)}
{\int_{\partial \Delta_Q(r)} M(x, z) \nu_1(dz)}
\,=\,c_1^2\,\frac{u_1(x)}{v_1(x)}\,\le\, c_1^4
\, \frac{u(x)}{v(x)}.
\end{eqnarray*}
Since $u_1, v_1$ are harmonic for $X^U$ and vanish
continuously on $\partial \Delta_Q(r) \cap \partial D$, by the
boundary Harnack principle
(Theorem 4.6 in \cite{KS1}), there exist $N$ and $c_2$ such that
$$
\frac{u_1(x)}{v_1(x)} \,\le\, c_2 \frac{u_1(y)}{v_1(y)}, \quad x,y \in
D \cap B(Q, \frac{r}N).
$$
Thus for every $ x,y \in
D \cap B(Q, \frac{r}N)$
$$
\frac{u(x)}{v(x)} \,\le\,c_1^2\,\frac{u_1(x)}{v_1(x)}\,\le\,c_2 \,c_1^2\,
 \frac{u_1(y)}{v_1(y)}\,\le\,c_2\, c_1^4 \,\frac{u(y)}{v(y)}.
$$
\qed

\vspace{.2in}
\begin{singlespace}

\end{singlespace}
\end{doublespace}
\end{document}